\documentclass[11pt,leqno]{article}
\usepackage{amssymb} 
\usepackage{epsfig}
\usepackage{mathrsfs} 
\usepackage{amsmath}

\newcommand{\bea}{\begin{eqnarray}}
\newcommand{\eea}{\end{eqnarray}}
\newcommand{\beas}{\begin{eqnarray*}}
\newcommand{\eeas}{\end{eqnarray*}}
\newcommand{\beq}{\begin{equation}}
\newcommand{\eeq}{\end{equation}}
\newcommand{\nn}{\nonumber}
\newcommand{\dsp}{\displaystyle}

\newcommand{\MM}{\mbox{$\mathcal M$}}

\newcommand{\GG}{\mbox{$\mathcal G$}}

\title{\textbf{
An ADI extrapolated Crank-Nicolson orthogonal spline collocation method
for nonlinear reaction-diffusion systems}}
\author{
Ryan I. Fernandes
\thanks{Department of Mathematics, The Petroleum Institute, P. O. Box 2533, Abu Dhabi, UAE.
(E-mail: rfernandes@pi.ac.ae). This author's research was supported by a research grant no. 12019 from the Petroleum Institute.}
\and
Graeme Fairweather
\thanks{Mathematical Reviews, American Mathematical Society, 416 Fourth Street, Ann Arbor,
MI 48103, U.S.A.
(Email: gxf@ams.org).}
}

\begin{document}
\maketitle

\begin{abstract}
An alternating direction implicit (ADI) orthogonal spline collocation (OSC) method is described for the approximate solution of a class of nonlinear reaction-diffusion systems. Its efficacy is demonstrated on the solution of
well-known examples of such systems, specifically the Brusselator,  Gray-Scott, Gierer-Meinhardt and Schnakenberg models,
and comparisons are made with other numerical techniques considered in the literature. The new ADI method is based on an extrapolated Crank-Nicolson OSC method and is algebraically linear.
It is efficient, requiring at each time level only $O({\cal N})$ operations where ${\cal N}$ is the number of unknowns. Moreover, it is shown to produce approximations which are of optimal global accuracy in various norms, and
to possess superconvergence properties.

\end{abstract}

{\small{\bf Key words}.
Nonlinear reaction-diffusion systems, alternating direction implicit method,
orthogonal spline collocation,
extrapolated Crank-Nicolson method, Brusselator, Gray-Scott, Gierer-Meinhardt, Schnakenberg models. }

\medskip
{\small {\bf AMS subject classifications.} 65M70 , 65M22, 92B05, 92-08}.

\maketitle

\setcounter{equation}{0}
\section{Introduction}
\label{sec:1}
In this paper, we consider  two-component reaction-diffusion problems
of the form:
\begin{eqnarray}
\label{prob1}
&& \frac{\partial {\bf u}}{\partial t} - {\cal D} \Delta {\bf u}
=  {\bf f}({\bf u}),
\quad (x,y,t) \in \Omega \times (0,T],\\
&&\frac{\partial {\bf u}}{\partial {{n}}} = {\bf 0}, \quad (x,y,t) \in \partial \Omega \times (0,T],\\
&&{\bf u}(x,y,0 ) = {\bf g}^{0} (x,y), \quad (x,y)
\in \overline{\Omega},
\label{prob4}
\end{eqnarray}
where $\Delta \equiv \partial^2/\partial x^2 + \partial^2/\partial y^2$  (the Laplace operator),
$\Omega$ is the rectangle $(a_1,b_1)\times(a_2,b_2)$ 
with boundary $\partial \Omega$,
$\overline{\Omega}=\Omega \cup \partial \Omega$,
${\bf u}=[u_1, u_2]^T$, 
${\cal D}=\mbox{diag}[D_1,D_2]^T$ 
is the diagonal matrix of diffusion constants, 
${\bf f}=[f_1,f_2]^T$ where the functions $f_1$ and $f_2$ specify the reaction kinetics of the system,
and $\partial/\partial n$ is the
outward normal derivative on the boundary $\partial \Omega$.
Problems of this type model phenomena arising in various fields such as  chemistry, physics, ecology and biology.  In particular, Turing \cite{Tu} proposed such a reaction-diffusion problem
to explain biological pattern formation. See \cite{KoMi} for a recent discussion of this pioneering work. Over the last several decades,  reaction-diffusion models have evolved in different fields and have come to be known by 
names such as the Brusselator model \cite{NiPr}, the Gray-Scott model \cite{GrSc1, GrSc2},
the Gierer-Meinhardt model \cite{GiMe}, and the Schnakenberg model \cite{Sc}. The derivation of these models is discussed in \cite{CrMa}, for example.

Alternating direction implicit (ADI) methods were first introduced in the context of finite
difference methods by Peaceman and Rachford \cite{PeRa} over fifty
years ago, and continue to be studied extensively for the solution of a variety of problems; see, for example,
\cite{BiFr,BrLy,Cu,DaChJa,GaXi,KnSu,LyBr,PoTo,Wat,Way,Ya} and references in these papers.
There has been much work done on the development, analysis and implementation of ADI orthogonal spline collocation (OSC) methods for the solution
of various scalar transient problems, an overview of which is given in \cite{FeBiFa}.
The purpose of this paper is to present an ADI extrapolated Crank-Nicolson (ECN)  OSC method (called simply the ADI method in the sequel) for the numerical solution of 
(\ref{prob1})--(\ref{prob4}). This ADI method is based on an ADI ECN OSC scheme formulated in \cite{BiFe1} for general nonlinear scalar parabolic problems.
Some of the salient features of the ADI method are the following.
\begin{itemize}

 \item
 It is a second-order accurate in time perturbation of a Crank-Nicolson OSC scheme and preserves the accuracy of 
that scheme.

\item
 By using extrapolation to linearize (\ref{prob1}), it is algebraically linear; that is, the algebraic systems to be solved at each step of the method are linear.

\item
Like all ADI methods, it reduces the solution of the multidimensional problem to the solution of sets of one-dimensional problems in each coordinate direction at each time step. With standard choices of bases for the spline space in the OSC spatial discretization, these problems involve almost block diagonal (ABD) linear systems  \cite{AmCa,FaGl}, which can be solved efficiently  using existing software \cite{DiFaKe1,DiFaKe2}.
Furthermore, since the differential operator in (\ref{prob1})
has constant coefficients, the ABD coefficient matrices in each coordinate direction are independent of time and can be decomposed only once for the entire problem.

\item Unlike  finite element Galerkin methods, OSC does not involve the approximation of integrals. Moreover, without post-processing, it provides superconvergent approximations at the nodes of
the partition of the spatial domain $\overline{\Omega}$ on which the spline space is defined. Unlike finite difference methods, it yields continuous approximations to the solution and its first derivatives of high-order accuracy throughout $\overline{\Omega}$.

\end{itemize}

A brief outline of the paper is as follows. Definitions and basic notation are given in section 2, followed by a description of the ADI scheme in section 3. In section 4, the results of numerical experiments with test problems from the literature involving the Brusselator, Gray-Scott, Gierer-Meinhardt and Schnakenberg models are presented, and comparisons are made with other published solution techniques. Concluding remarks are given in section 5.

\setcounter{equation}{0}
\section{Preliminaries}
\label{sec:2}
Let
$ 
\{x_i\}_{i=0}^{N_x}$ and $
\{y_j\}_{j=0}^{N_y}$
be partitions (in general, non-uniform) of
$[a_1,b_1]$ and $[a_2,b_2]$, respectively, 
such that
\[
a_1 = x_0<x_1<\cdots<x_{ _{ N_x-1}}<x_{ _{N_x}}=b_1,\hspace{.2in}
a_2 = y_0<y_1<\cdots<y_{ _{N_y-1}}<y_{ _{N_y}}=b_2.
\]
Let $\MM_x$ and $\MM_y$
be the spaces of piecewise
polynomials of degree $\leq r$, $r \geq 3$, defined by
\beas
\MM_x &=& \{v \in C^1[a_1,b_1]: v|_{[x_{i-1},x_i]} \in P_r,\ 1 \le i \le N_x\},
\\
\MM_y &=& \{v \in C^1[a_2,b_2]: v|_{[y_{j-1},y_j]} \in P_r,\ 1 \le j \le N_y\},
\eeas
where $P_r$ denotes the set of polynomials of degree $\leq r$.
Let
$
\MM = \MM_x \otimes \MM_y,
$
the set of all functions that are finite linear combinations of
products of the form $v^x(x)v^y(y)$ where $v^x \in \MM_x$ and $v^y \in \MM_y$.
Let $\{ \lambda_k \}_{k=1}^{r-1}$ 
be the nodes 
of the $(r-1)$-point
Gauss quadrature rule on $[0,1]$ 
and define
$\GG_x=\{\xi^x_{i,k}\}_{i=1,k=1}^{\;N_x, r-1}$,
$\GG_y=\{\xi^y_{j,l}\}_{j=1,l=1}^{\;N_y, r-1}$, where
\[
\xi^x_{i,k}=x_{i-1}+h^x_i\lambda_k,\hspace{.2in}
\xi^y_{j,l}=y_{j-1}+h^y_j\lambda_l, \quad h_i^x=x_i-x_{i-1}, \quad h_j^y=y_j-y_{j-1}.
\]
Then
\[\GG = \{ \xi = (\xi^x, \xi^y ): \xi^x \in \GG_x, \xi^y \in \GG_y\}
\]
is the set of 
collocation points in $\Omega$.
Let $\{t_n\}_{n=0}^{N_t}$ be a partition of $[0,T]$ such that $t_n=n \tau,$
where $\tau=T/N_t,$ and, for $n=1,\ldots,N_t-1$, let $t_{n+1/2}=(t_n+t_{n+1})/2$.
Throughout this paper, for any function $\phi$ of $t$, $\phi^n(\cdot)=\phi(\cdot,t_n)$.

\section{The ADI Extrapolated Crank-Nicolson OSC Method}
\label{sec:3}
\setcounter{equation}{0}
In this section, we describe the ADI  scheme. As we shall see, at each time level, the approximate solution is a
piecewise polynomial of degree $\leq r$ in one variable only along lines through the
collocation points. At the final step (or at any user-specified intermediate step), the approximation
is converted to a two dimensional piecewise polynomial on $\overline{\Omega}$.

We compute the OSC approximation  ${\bf u}_h =[u_{h,1},u_{h,2}]^T \in \MM\times\MM$ to the
solution ${\bf u}$ of (\ref{prob1})--(\ref{prob4}) at the final time $T$ in the following way.

\vspace{.2in}
\noindent
{\bf Step 1a.}
First, we calculate starting values by interpolation at Gauss points.

\vspace{.1in}
\noindent
$\bullet$ \hspace{.1in} For each $\xi^x \in \GG_x$, we choose ${\bf u}_h^0(\xi^x,\cdot) \in \MM_y\times\MM_y$
(along {\it vertical} lines through the collocation points) such that
\beq
\label{eq:iv1}
{\bf u}_h^0 (\xi^x,\widehat{\xi^y}) = {\bf g}^0(\xi^x,\widehat{\xi^y}), \quad \widehat{\xi^y} \in \GG_y \cup \{a_2, b_2 \}. 
\eeq
$\bullet$  \hspace{.1in} For each $\xi^y \in \GG_y$, we choose $\widehat{\bf u}_h^0 (\cdot,\xi^y) \in \MM_x\times\MM_x$
(along {\it horizontal} lines through the collocation points) such that
\beq
\label{eq:iv2}
\widehat{\bf u}_h^0 (\widehat{\xi^x},\xi^y) = {\bf g}^0(\widehat{\xi^x},\xi^y), \quad \widehat{\xi^x} \in \GG_x \cup \{a_1, b_1 \}. 
\eeq

\vspace{.2in}
\noindent
{\bf Step 1b.}
We also require a second-order accurate
approximation to the exact solution ${\bf u}(x,y,t_{1/2})$, which we define in the following way. Using Taylor's theorem, (\ref{prob1})
with $t=0$ and (\ref{prob4}), we have,
\[
{\bf u}(x,y,t_{1/2})={\bf g}^0(x,y) + \frac{\tau}{2} [{\bf f}({\bf g}^0(x,y))+{\cal D} \Delta{\bf g}^0(x,y)] + O(\tau^2),
\quad (x,y) \in \Omega.
\]
Then, for each $\xi^x \in \GG_x$, we choose
$\widetilde{\bf u}_h^{1/2}(\xi^x,\cdot) \in \MM_y\times\MM_y$ (along {\it vertical} lines through the collocation points) such that
\beq
\label{eq:iv3}
\left\{
\begin{array}{l}
\widetilde{\bf u}_h^{1/2} (\xi^x,\xi^y) = {\bf g}^0(\xi^x,\xi^y)+\dsp{\frac{\tau}{2} \left[ {\bf f}({\bf u}_h^0(\xi^x,\xi^y))
+{\cal D} \left ( \frac{\partial^2 \widehat{\bf u}_h^0}{\partial x^2}
+\frac{\partial^2 {\bf u}_h^0}{\partial y^2} \right)(\xi^x,\xi^y)
\right]}, \quad \xi^y \in \GG_y, \\ \\
\dsp{\frac{\partial \widetilde{\bf u}_h^{1/2}}{\partial y} (\xi^x,\alpha)=0, \quad \alpha=a_2, b_2.} 
\end{array}
\right.
\eeq

With standard choices of bases for $\MM_x$ and $\MM_y$, the interpolation processes
(\ref{eq:iv1}), (\ref{eq:iv2}) and (\ref{eq:iv3}) involve
ABD 
linear systems of orders $N_y(r-1)+2$, $N_x(r-1)+2$ and $N_y(r-1)+2$, respectively, each of which can be solved at a cost of $O(N_x)$ or $O(N_y )$; see \cite{DiFaKe1}.

Next, we advance in time; for each $n=0,1,\ldots,N_t-1$, we perform {\bf Step 2} followed by {\bf Step 3} as follows:

\vspace{.2in}
\noindent
{\bf Step 2.}
For each $\xi^y \in \GG_y$, determine ${\bf u}_h^{n+1/2}(\cdot,\xi^y) \in \MM_x\times\MM_x$ such that
\beq
\left\{
\begin{array}{l}
\dsp{\left [\frac{{\bf u}^{n+1/2}_h-{\bf u}^n_h}{\tau/2}-{\cal D}\left (\frac{\partial^2{\bf u}_h^{n+1/2}}{\partial x^2}
+\frac{\partial^2{\bf u}_h^{n}}{\partial y^2}\right )\right](\xi^x,\xi^y) = {\bf f}(\widetilde{\bf u}_h^{n+1/2}(\xi^x,\xi^y)),
\quad\xi^x \in \GG_x}, \\ \nn \\ \nn
\dsp{\frac{\partial {\bf u}_h^{n+1/2}}{\partial x}} (\alpha,\xi^y)=0, \quad \alpha=a_1, b_1, 
\end{array}
\right.
\eeq
where $\widetilde{\bf u}_h^{n+1/2}(\xi^x,\cdot) \in \MM_y\times\MM_y$ is given by
\beq
\label{eq:extrap}
\widetilde{\bf u}_h^{n+1/2}(\xi^x,\cdot)=[ 3 {\bf u}_h^n (\xi^x,\cdot)-{\bf u}_h^{n-1} (\xi^x,\cdot)]/2,
\quad n=1,\ldots,N_t-1,
\eeq
and ${\bf u}_h^0(\xi^x,\cdot)$, $\widetilde{\bf u}_h^{1/2}(\xi^x,\cdot)$ are given by (\ref{eq:iv1}),
(\ref{eq:iv3}), respectively.

\vspace{.1in}
\noindent
This step comprises a set of independent OSC two-point boundary value problems (TPBVPs) along {\it horizontal} lines through the
collocation points ${\cal G}_y$, each of which gives rise to an ABD linear system of order $N_y(r-1)+2$.

\vspace{.2in}
\noindent
{\bf Step 3.} For each $\xi^x \in \GG_x$, determine ${\bf u}_h^{n+1}(\xi^x,\cdot) \in \MM_y\times\MM_y$  such that
\beq
\left\{
\begin{array}{l}
\dsp{\left [\frac{{\bf u}_h^{n+1}-{\bf u}_h^{n+1/2}}{\tau/2}-{\cal D}\left ( \frac{\partial^2{\bf u}_h^{n+1/2}}{\partial x^2}
+\frac{\partial^2{\bf u}_h^{n+1}}{\partial y^2}\right )\right ] (\xi^x,\xi^y) = {\bf f}(\widetilde{\bf u}_h^{n+1/2}(\xi^x,\xi^y)),
\quad \xi^y \in \GG_y,} \\ \nn \\ \nn
\dsp{\frac{\partial {\bf u}_h^{n+1}}{\partial y}}
(\xi^x,\alpha)=0, \quad \alpha=a_2, b_2,
\end{array}
\right.
\eeq
with $\widetilde{\bf u}_h^{n+1/2}(\xi^x,\cdot)$ as in (\ref{eq:extrap}).

\vspace{.1in}
\noindent
This step comprises a set of independent OSC TPBVPs along {\it vertical} lines through the
collocation points ${\cal G}_x$, each of which gives rise to an ABD linear system of order $N_x(r-1)+2$.

\vspace{.2in}
Once these two steps are completed, we proceed to the final phase of the algorithm,
the purpose of which is to convert the one dimensional approximations along vertical lines through the
collocation points ${\cal G}_x$,
${\bf u}_h^{N_t}(\xi^x,\cdot) \in \MM_y\times\MM_y,$  $\xi^x \in \GG_x,$ to the desired two dimensional approximate solution ${\bf u}_h (\cdot,\cdot,T) \in \MM\times\MM$
on $\overline{\Omega}$. To this end, we observe that ${\bf {u}}_h (\xi^x,\cdot,T) \in \MM_y \times \MM_y$ for $\xi^x \in \GG_x,$ and, for fixed $\widehat{\xi}^y \in \GG_y \cup \{a_2,b_2\},$ we introduce ${\bf \overline{u}}_h (\cdot,\widehat{\xi}^y,T) \in \MM_x \times \MM_x$. Then, to obtain the desired approximate solution, we perform the following steps.

\vspace{.2in}
\noindent
{\bf Step 4a.} For each $\widehat{\xi}^y \in \GG_y\cup\{a_2,b_2\}$, 
determine ${\bf \overline{u}}_h (\cdot,\widehat{\xi}^y,T) \in \MM_x \times \MM_x$ along \textit{horizontal} lines through $\GG_y\cup\{a_2,b_2\}$ such that
\[
\left\{
\begin{array}{l}
{\bf \overline{u}}_h(\xi^x,\widehat{\xi}^y,T) = {\bf u}_h^{N_t}(\xi^x,\widehat{\xi}^y), \quad \xi^x \in \GG_x, 
\\ \\
\dsp{\frac{\partial {\bf \overline{u}}_h}{\partial x}}
(\alpha,\widehat{\xi^y},T) = 0, \quad \alpha=a_1, b_1. 
\end{array}
\right.
\]

\vspace{.2in}
\noindent
{\bf Step 4b.}
Next, determine ${\bf u}_h(\alpha,\cdot,T) \in \MM_y \times \MM_y$, where $\alpha =a_1,b_1$, along \textit{vertical} lines through $\{a_1,b_1\}$ such that
\[
{\bf u}_h(\alpha,\widehat{\xi}^y,T) = {\bf \overline{u}}_h(\alpha,\widehat{\xi}^y,T), \quad 
\widehat{\xi}^y \in \GG_y\cup\{a_2,b_2\}.
\]

\vspace{.2in}
\noindent
{\bf Step 4c.}
Lastly, determine ${\bf u}_h(\cdot,\cdot,T) \in \MM \times \MM$, along \textit{horizontal} lines through $y \in [a_2, b_2]$ such that
\beq
\label{eq:4c}
{\bf u}_h(\xi^x,y,T) = {\bf u}_h^{N_t}(\xi^x,y),
\quad  \xi^x \in \GG_x, 
\eeq
where the boundary values for (\ref{eq:4c}) are obtained in Step 4b. 

\vspace{.1in}
\noindent
Steps 4a and 4c  involve the solution of ABD systems of order $N_x(r-1)+2$ whereas Step 4b requires the solution of two ABD systems of order $N_y(r-1)+2$.
While we have assumed that the approximate solution is required only at the final value $T$, Steps 4a, b, c can be performed at any intermediate time level, $t_n$, $1\le n \le N_t-1$.

It can be shown that the total cost of computing the approximate solution ${\bf u}_h$ using the ADI scheme is $O(r^4 N_x N_y N_t)$, where $r$ is the degree of the piecewise polynomials and $N_x$, $N_y$ and $N_t$ are the numbers of subdivisions  in the $x$, $y$ and $t$ directions, respectively. Details of the implementation of the method are similar to those of \cite[Section 2.2]{BiFe1} and are
omitted.

\section{Numerical Experiments}
\label{sec:4}
\setcounter{equation}{0}

In this section, we consider commonly occurring  reaction-diffusion models of the form (\ref{prob1})--(\ref{prob4}) and compare results obtained by applying the ADI scheme
with results appearing in the literature. We also demonstrate the global accuracy and superconvergence properties of the ADI scheme. In all of the test problems, we consider uniform partitions in the $x$ and $y$ directions with
$N_x=N_y=N$.

\subsection{Brusselator Model}
\label{bm}
The system of partial differential equations associated with the Brusselator model \cite{NiPr} comprises (\ref{prob1}) with
\[
{\bf f}({\bf u}) =[B + u_1^2 u_2 - (A + 1)u_1,A u_1 - u_1^2 u_2]^T,
\]
where 
$A,B$ are constants. 
Various numerical methods have been proposed for its solution. Twizell et al., \cite{TwGuCa} formulated an algebraically linear finite difference method which is second-order accurate in space and time. Adomian \cite{Ad} proposed a decomposition method
which was subsequently corrected and extended in \cite{Wa}. Verwer et al.,
\cite{VeHuSo} discretized in space using a basic finite difference approximation and employed a conditionally stable, explicit second-order Runge-Kutta-Chebyshev method for the time-stepping.  Ang \cite{An} formulated a dual-reciprocity boundary element method which is claimed to be applicable to problems in domains of arbitrary shape.  More recently,
in \cite{SiAlHa}, a method
involving a spatial discretization based on radial basis functions and a linearized Crank-Nicolson-type method for the time-stepping is formulated, while in \cite{MiJi1,MiJi2}
differential quadrature is used in space, and in time a fourth-order  accurate Runge-Kutta method of Pike and Roe \cite{PiRo}.

\vspace{.2in}
\noindent
{\bf Example 1}. To demonstrate the optimal accuracy of the ADI  scheme in various norms,
we consider a test problem in which
\[
\Omega=(0,1)\times (0,1), \quad
 D_1=D_2=1, \quad A=1,\quad B=0.5, \quad T=1,
 \]
and the exact solution of (\ref{prob1})--(\ref{prob4}) has components
 \beq
 \label{exsol}
u_1(x,y,t)= \cos(t) \cos(2 \pi x) \cos(\pi y), \qquad
u_2(x,y,t)= \cos(t) \cos(\pi x) \cos(2 \pi y).
\eeq
Clearly ${\bf g}^0(x,y) = [u_1(x,y,0),u_2(x,y,0)]^T$, and
we construct the reaction kinetic functions in (\ref{prob1})
to be
\[
\widetilde{\bf f}({\bf w})=\widehat{\bf f}(x,y,t)+{\bf f}({\bf w}),
\]
for any  ${\bf w} =[w_1,w_2]^T$,
where
\[
\widehat{\bf f}(x,y,t) =\frac{\partial {\bf u}}{\partial t} - {\cal D} \Delta {\bf u} -
{\bf f}({\bf u}).
\]

In the OSC discretization, we consider splines of degree $r$, where
$r = 3,4,5,$  and choose the time step  $\tau=h^{(r+1-k)/2}, k=0,1,$
since it is expected that
the method is $O(\tau^2+h^{r+1-k})$ accurate in the $H^k(\Omega)$ norm. 
For various values of $N$ and $T=1$, Tables 1 and 2 show, respectively, the $L^2(\Omega)$ and $H^1(\Omega)$ norms of the errors in $u_{h,1}$
and $u_{h,2}$, denoted by $e_{h,1}$ and $e_{h,2}$.
These norms of the errors are calculated using $(r+2)$-point Gauss quadrature so that the quadrature error does not affect the accuracy of the scheme. The tables also show the experimental convergence rate of the scheme, denoted by Rate, calculated using the formula
\begin{eqnarray*}
\mbox{Rate} =
\frac{\log(\mbox{error}_{h_1}/\mbox{error}_{h_2})}{\log(h_1/h_2)},
\end{eqnarray*}
where $\mbox{error}_{h}=\sqrt{\|e_{{h,1}}\|_{H^{k}(\Omega)}^2+\|e_{{h,2}}\|_{H^{k}(\Omega)}^2}$. %
The results in Tables 1 and 2 confirm the expected convergence rates of the  ADI scheme in the $H^k(\Omega)$ norm, $k=0,1$.
 \begin{table}[hbt]
\begin{center}
{\scriptsize
\begin{tabular}{|c|c|c|c|c|}
\hline $r$ & $N$  & $\|e_{h,1}\|_{L^2(\Omega)}$ & $\|e_{h,2}\|_{L^2(\Omega)}$ & Rate \\
\hline
3 & 10    &  0.683--04& 0.706--04&   \\ 
  & 15    &  0.134--04& 0.139--04& 4.014  \\ 
  & 20    &  0.424--05& 0.438--05& 4.007  \\  \hline \hline
4 & 9    &  0.139--04& 0.144--04&   \\ 
  & 16    &  0.785--06& 0.816--06& 4.993  \\ 
  & 25    &  0.845--07& 0.878--07& 4.996  \\  \hline \hline
5 & 10    &  0.818--06& 0.850--06&   \\ 
  & 15    &  0.718--07& 0.747--07& 6.000  \\ 
  & 20    &  0.128--07& 0.133--07& 6.000  \\  \hline
\end{tabular}
}
\end{center}
\begin{center}
\textbf{Table 1.} Example 1: $L^2(\Omega)$ errors  at $T=1$
\end{center}
\end{table}
 \begin{table}[hbt]
\begin{center}
{\scriptsize
\begin{tabular}{|c|c|c|c|c|}
\hline $r$ & $N$  & $\|e_{h,1}\|_{H^{1}(\Omega)}$ & $\|e_{h,2}\|_{H^{1}(\Omega)}$  & Rate \\
\hline
3 & 9     &  0.835--02& 0.864--02&   \\ 
  & 16    &  0.151--02& 0.156--02& 2.975  \\ 
  & 25    &  0.399--03& 0.411--03& 2.986  \\  \hline \hline
4 & 10    &  0.586--03& 0.609--03&   \\ 
  & 15    &  0.116--03& 0.120--03& 3.999  \\ 
  & 20    &  0.366--04& 0.381--04& 3.999  \\  \hline \hline
5 & 9     &  0.983--04& 0.102--03&   \\ 
  & 16    &  0.553--05& 0.576--05& 5.000  \\ 
  & 25    &  0.594--06& 0.618--06& 5.000  \\  \hline
\end{tabular}
}
\end{center}
\begin{center}
\textbf{Table 2.} Example 1: $H^1(\Omega)$ errors at $T=1$
\end{center}
\end{table}
 \begin{table}[hbt]
\begin{center}
{\scriptsize
\begin{tabular}{|c|c|c|c|c|}
\hline $r$ & $N$  & $\|e_{h,1}\|_{L^\infty(\Omega)}$  & $\|e_{h,2}\|_{L^\infty(\Omega)}$& Rate \\
\hline
3 & 10    &  0.239--03& 0.244--03&   \\ 
  & 15    &  0.472--04& 0.478--04& 4.019  \\ 
  & 20    &  0.149--04& 0.152--04& 3.985  \\  \hline \hline
4 &  9    &  0.282--04& 0.297--04&   \\ 
  & 16    &  0.162--05& 0.168--05& 4.999  \\ 
  & 25    &  0.174--06& 0.180--06& 4.997  \\  \hline \hline
5 & 10    &  0.172--05& 0.177--05&   \\ 
  & 15    &  0.152--06& 0.156--06& 5.997  \\ 
  & 20    &  0.270--07& 0.279--07& 5.983  \\  \hline
\end{tabular}
}
\end{center}
\begin{center}
\textbf{Table 3.} Example 1: $L^\infty(\Omega)$ errors at $T=1$
\end{center}
\end{table}

With $\tau=h^{(r+1)/2},$ we also compute the $L^{\infty}(\Omega)$ norm of the errors $e_{{h,1}}$ and $e_{{h,2}}$ at $T=1$ and the corresponding convergence rates. These $L^{\infty}(\Omega)$ norms are calculated by taking the maximum absolute values of $e_{{h,1}}$ and $e_{{h,2}}$ at $10 \times 10$ equally spaced points in each cell $[x_{i-1},x_i] \times [y_{j-1},y_j]$, $1 \leq i,j \leq N$. In Table 3, the $L^{\infty}(\Omega)$ convergence rate is seen to be optimal, that is, $O(\tau^2+h^{r+1})$. Note that, when computing Rate in this case, we take $\mbox{error}_h=\mbox{max}\{ \|e_{{h,1}}\|_{L^{\infty}},\| e_{{h,2}}\|_{L^{\infty}}\}$.

Invariably, OSC schemes exhibit superconvergence phenomena; specifically, one obtains an accuracy of $O(\tau^2+h^{2r-2})$ 
in the solution and the first spatial partial derivatives at the nodal points $\{(x_i,y_j)\}_{i,j=0}^{N}$.
To examine the superconvergence of the ADI  scheme, we choose $\tau=h^{r-1}.$ From the results presented in Table 4,
it is evident that the ADI scheme possesses the expected properties.

 \begin{table}[hbt]
\begin{center}
{\scriptsize
\begin{tabular}{|c|c|c|c|c|c|c|c|c|c|c|}
\hline
&    & \multicolumn{2}{c|}{Maximum nodal errors}  &   & \multicolumn{2}{c|}{Maximum nodal errors} &   & \multicolumn{2}{c|}{Maximum nodal errors} &  \\ \cline{3-4} \cline{6-7} \cline{9-10}
$r$ & $N$  & $e_{h,1}$ & $e_{h,2}$ & Rate & $\partial e_{h,1}/\partial x$ & $\partial e_{h,2}/\partial x$ & Rate & $\partial e_{h,1}/\partial y$ & $\partial e_{h,2}/\partial y$ & Rate\\
\hline
3 & 10 & 0.239--03& 0.244--03&       & 0.951--03& 0.752--03&      & 0.729--03& 0.983--03&  \\
  & 15 & 0.472--04& 0.478--04& 4.019 & 0.196--03& 0.146--03& 3.897& 0.143--03& 0.203--03& 3.896\\
  & 20 & 0.149--04& 0.152--04& 3.985 & 0.623--04& 0.468--04& 3.982& 0.454--04& 0.644--04& 3.985 \\  \hline \hline
4 & 10 & 0.160--05& 0.165--05&       & 0.969--05& 0.516--05&      & 0.489--05& 0.101--04&  \\
  & 15 & 0.141--06& 0.144--06& 6.024 & 0.887--06& 0.443--06& 5.895& 0.427--06& 0.928--06& 5.894\\
  & 20 & 0.250--07& 0.258--07& 5.966 & 0.159--06& 0.806--07& 5.979& 0.764--07& 0.166--06& 5.984 \\  \hline \hline
5 & 10 & 0.167--07& 0.172--07&       & 0.986--07& 0.542--07&      & 0.513--07& 0.103--06&  \\
  & 15 & 0.650--09& 0.665--09& 8.026 & 0.401--08& 0.206--08& 7.896& 0.199--08& 0.417--08& 7.897\\
  & 20 & 0.666--10& 0.685--10& 7.901 & 0.413--09& 0.215--09& 7.906& 0.205--09& 0.428--09& 7.913 \\
\hline
 \end{tabular}
 }
\end{center}
\begin{center}
\textbf{Table 4.} Example 1: Maximum nodal errors at $T=1$
\end{center}
\end{table}

\vspace{.2in}
\noindent
{\bf Example 2.}
In this test problem from \cite{TwGuCa}, which has no known closed-form solution, 
we have
\[
\Omega=(0,1)\times (0,1), \quad
 D_1=D_2=0.002, \quad A=1,\quad B=2, \quad T=5,
 \]
 and
\[
{\bf g}^0(x,y) = [2+0.25 y,1+0.8 x]^T.
\]
Figure 1 shows the initial ($t=0$) and final ($t=T=5$) graphs of $u_{h,1}$ 
and $u_{h,2}$
obtained using the ADI scheme with
$N=10$, (that is, $h=0.1$), $r = 3$ and $\tau=h^{2}$.
The results for $T=5$ (and for $T=10$ which are similar, but not presented here) are in agreement with the theory of \cite{TwGuCa} which predicts that the solution ${\bf u}$ converges to the fixed point $(B, A/B)$ when
\beq
\label{eq:constraint}
1-A+B^2 \ge 0.
\eeq

The graphs in Figure 1 should compare with Figures 7 and 8 in \cite{TwGuCa} and Figures 4 and 5 in \cite{SiAlHa}. However, the figures in these papers are in error since they present graphs of approximations to $u_1(y,x)$ and $u_2(y,x)$ instead of $u_1(x,y)$ and $u_2(x,y)$, respectively.\footnote{E. H. Twizell, {\em Personal communication}, July 2011.}
\footnote{ Siraj-ul-Islam, {\em Personal communication}, July 2011.}

\begin{center}
\resizebox{1\linewidth}{0.7\linewidth}{\includegraphics{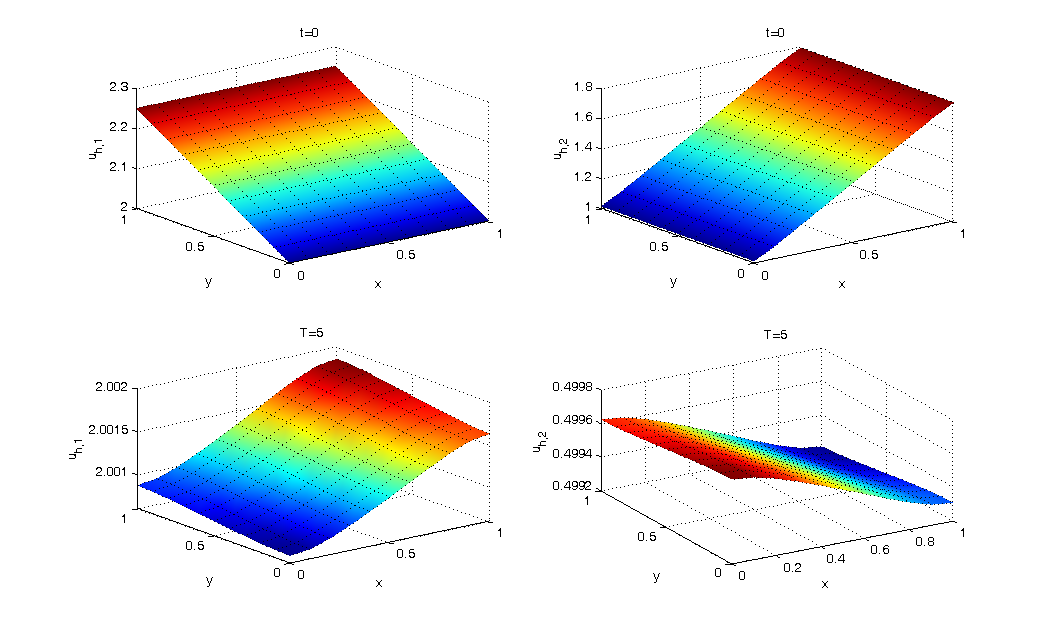}}
\end{center}
\begin{center}
{\bf Figure 1.}  Example 2: The initial ($t=0$) and final ($T=5$) approximate solutions, $u_{h,1}$ and $u_{h,2}$, 
\end{center}

\vspace{.2in}
 \noindent
 {\bf Example 3.} This test problem, presented in \cite{Ad} with no numerical results,
is similar to Example 2 but with $A=3.4$, $B=1$, $T=1$.
With $N=10$, $r = 3$ and $\tau=h^{2}$ as in Example 2, the ADI approximate solutions $u_{h,1}$ and $u_{h,2}$ at $T=1$ are shown in Figure 2. It should be noted that they
do not exhibit the oscillatory behavior present in \cite[Figures 9, 10]{TwGuCa}. This test problem is also solved successfully in \cite{SiAlHa}, and if Figure 8 of
\cite{SiAlHa} were
correctly drawn, it would be similar to Figure 2 of the present paper. Moreover, the problem is solved in \cite{Wa}
but little evidence is presented to demonstrate the efficacy of the method described therein.
\begin{center}
\resizebox{1\linewidth}{0.4\linewidth}{\includegraphics{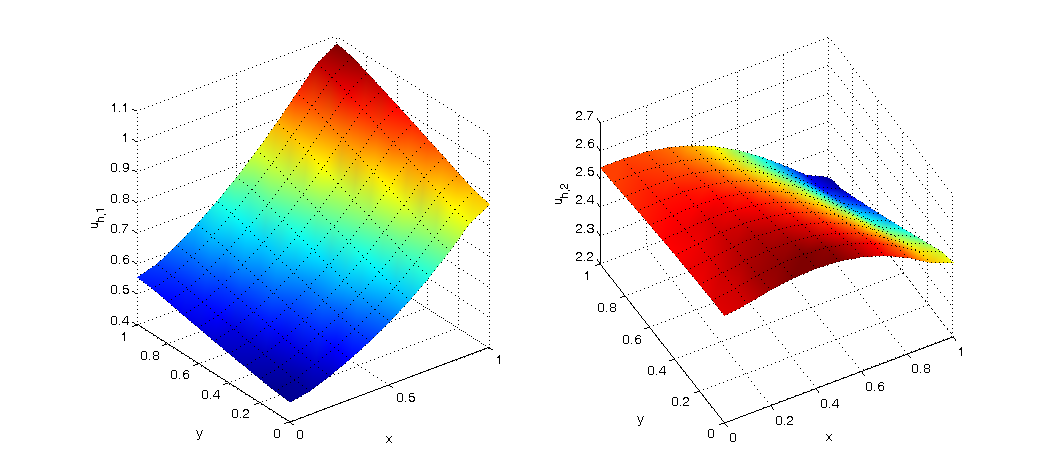}}
\end{center}
\begin{center}
{\bf Figure 2.} Example 3: The approximate solutions $u_{h,1}$ and $u_{h,2}$ at $T=1$
\end{center}

With this choice of $A$ and $B$, condition (\ref{eq:constraint}) is violated and, not surprisingly,  convergence to a fixed point is not observed in numerical experiments.
For example,  graphs of the approximate solutions $u_{h,1}$ and $u_{h,2}$ at $T=20$ and $40$ are given in Figure 3; cf., \cite[Figure 9]{SiAlHa}.
To better understand the nature of the solutions as $t$ increases,
we plot in Figure 4 the profiles of $u_{h,1}$ and $u_{h,2}$ at some points in $\Omega$ as $t$ ranges from $0$ to $40$. From this figure, it is obvious that the solutions
are oscillatory; see also \cite[Figure 9]{SiAlHa}.

\begin{center}
\resizebox{0.8\linewidth}{!}{\includegraphics{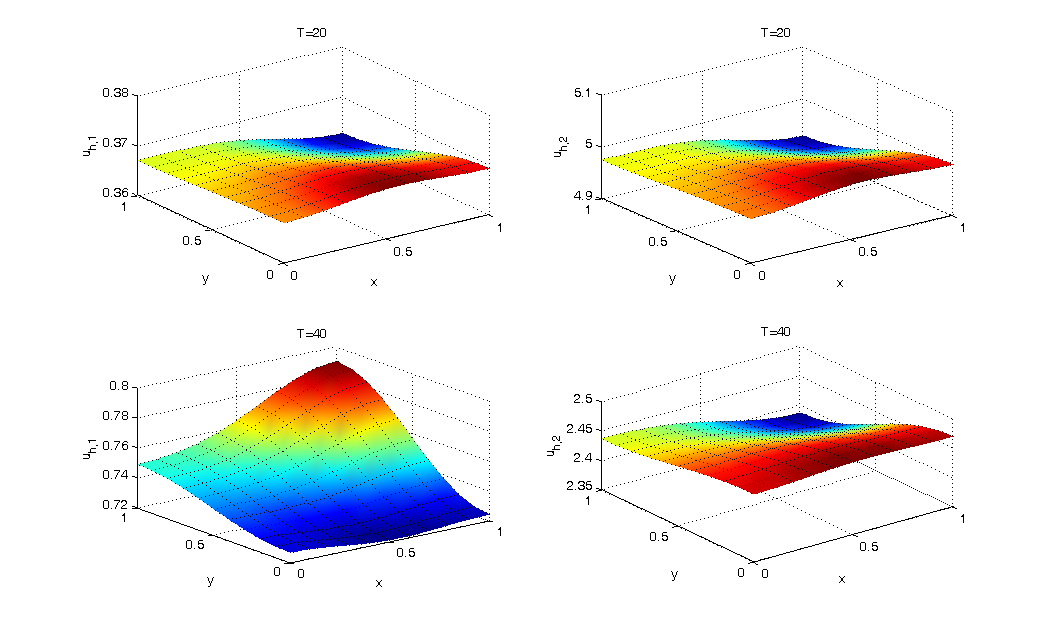}} 

{\bf Figure 3.} Example 3: The approximate solutions $u_{h,1}$ and $u_{h,2}$ at $T=20, 40$
\end{center}
\begin{center}
\resizebox{0.8\linewidth}{!}{\includegraphics{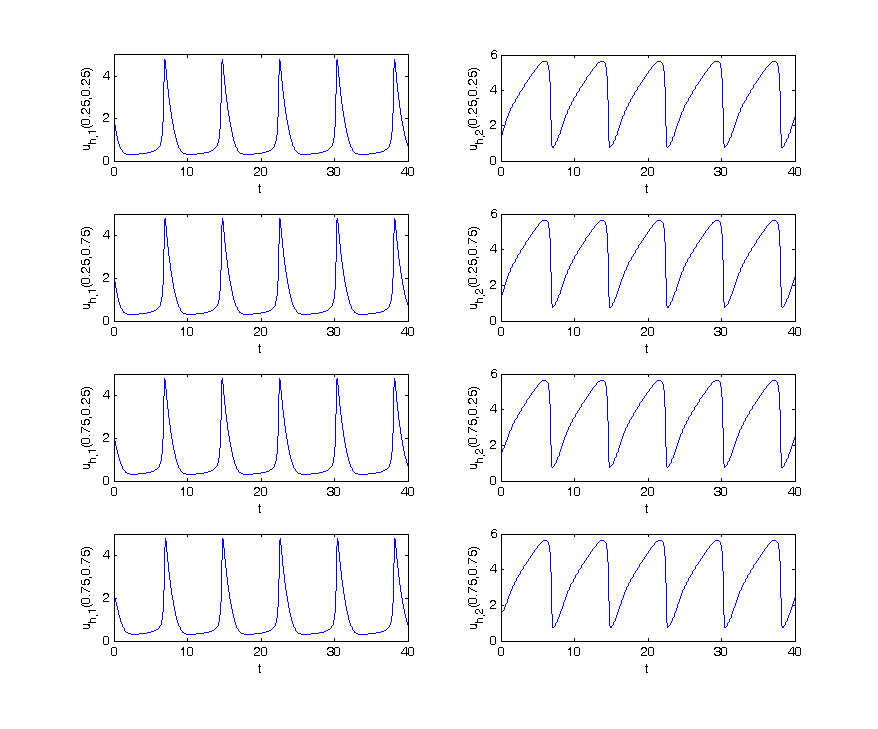}}
\end{center}
\begin{center}
{\bf Figure 4.} Example 3: The approximate solutions $u_{h,1}$ and $u_{h,2}$ against $t$ at $(0.25,0.25),$ $(0.25,0.75),$ $(0.75,0.25),$ $(0.75,0.75)$ 
\end{center}

Additional test problems appear in the literature.  Ang \cite{An} considered a problem similar to Example 2 but with
\[
{\bf g}^0(x,y) = [\frac{1}{2}x^2 -\frac{1}{3}x^3, \frac{1}{2}y^2 -\frac{1}{3}y^3]^T.
\]
However, few numerical results are presented in that paper. Verwer et al., \cite{VeHuSo} considered a test problem from \cite{HaNoWa}
which has the same selection of parameter values as in Example 3 but with
\[
{\bf g}^0(x,y) = [2+0.25 y,1+0.8 x]^T.
\]
They presented graphs of their approximation to $u_{2}$ at various time levels.
In \cite{MiJi1,MiJi2}, which are essentially identical papers, the authors claim to consider the same test problem but choose $A=1$ and $B=3.4$ instead of $A=3.4$ and $B=1$.
They also claim that graphs of their approximation to $u_1$ are similar to those in \cite{VeHuSo} when, in fact, \cite{VeHuSo} contains graphs of approximations to $u_2$ only.
Ang's test problem and Example 2 of the present paper are also considered in \cite{MiJi1,MiJi2}.  However, all of the graphs presented in \cite{MiJi1,MiJi2} are incorrect, since,
in each case, approximations to $u_1(y,x)$ and $u_2(y,x)$ instead of $u_1(x,y)$ and $u_2(x,y)$, respectively, are plotted.

\subsection{Gray-Scott Model}
The Gray-Scott model \cite{GrSc1,GrSc2} is given by
(\ref{prob1})--(\ref{prob4}) with
\[
{\bf f}({\bf u}) = [F(1-u_1)- u_1^2 u_2, u_1^2 u_2-(F+k)u_2]^T,
\]
where  $F$ and $k$ are constants.

\vspace{.2in}
\noindent
{\bf Example 4.}
This test problem is Example 2 in \cite{ZhWoZh} in which $\Omega = (-1,1)\times(-1,1)$, $T=1$, and the exact solution is chosen to be
\[
u_1(x,y,t)= \cos(2t) \cos(2 \pi x) \cos(\pi y), \qquad
u_2(x,y,t)= \cos(2t) \cos(\pi x) \cos(2 \pi y).
\]
Two parameter choices are considered:
\beq
\label{eq:case2}
D_1=D_2=0.001, \quad F=1, \quad k=0,
\eeq
and
\beq
\label{eq:case1}
 D_1=D_2=1, \quad F=1, \quad k=0.
\eeq
Zhang et al., \cite{ZhWoZh} used a finite element Galerkin method for the spatial discretization with spaces of $C^0$ piecewise linear or $C^0$ piecewise quadratic elements defined on uniform triangulations of $\overline{\Omega}$, and a time-stepping procedure based on a linearized second-order backward differentiation formula (BDF).  These methods yield approximations which are second- and third-order in space, respectively, and second-order in time. At each time step, the linear algebraic system is solved using an iterative method based on GMRES.

Results obtained for (\ref{eq:case2}) using the ADI method with a choice of $\tau$ consistent with the expected spatial accuracy are presented in Tables 5--7.  These tables demonstrate the optimal convergence rates in the $H^k(\Omega), k=0,1,$ and $L^{\infty}(\Omega)$ norms, respectively. Note that  $N$, the number of space intervals, is chosen so that  $N_T$, the number of time steps, is an integer. Superconvergence at the nodes is exhibited in Table 8 for (\ref{eq:case1}) which is the parameter choice used in \cite{ZhWoZh} to examine rates of convergence.
 \begin{table}[hbt]
\begin{center}
{\scriptsize
\begin{tabular}{|c|c|c|c|c|}
\hline $r$ & $N$  & $\|e_{h,1}\|_{L^2(\Omega)}$ & $\|e_{h,2}\|_{L^2(\Omega)}$ & Rate \\
\hline
3 & 20    &  0.708--04& 0.647--04&   \\ 
  & 26    &  0.252--04& 0.231--04& 3.928  \\ 
  & 32    &  0.111--04& 0.102--04& 3.963  \\  \hline 
4 &  8    &  0.569--03& 0.458--03&   \\ 
  & 18    &  0.984--05& 0.788--05& 5.005  \\ 
  & 32    &  0.567--06& 0.459--06& 4.953  \\  \hline \hline
5 & 20    &  0.575--06& 0.456--06&   \\ 
  & 26    &  0.119--06& 0.945--07& 6.000  \\ 
  & 32    &  0.343--07& 0.272--07& 6.000  \\  \hline
\end{tabular}
}
\end{center}
\begin{center}
\textbf{Table 5.} Example 4: $L^2(\Omega)$ errors at $T=1$ for (\ref{eq:case2})
\end{center}
\end{table}
 \begin{table}[hbt]
\begin{center}
{\scriptsize
\begin{tabular}{|c|c|c|c|c|}
\hline $r$ & $N$  & $\|e_{h,1}\|_{H^{1}(\Omega)}$ & $\|e_{h,2}\|_{H^{1}(\Omega)}$& Rate \\
\hline
3 & 8     &  0.820--01& 0.722--01&   \\ 
  & 18    &  0.756--02& 0.676--02& 2.932  \\ 
  & 32    &  0.137--02& 0.123--02& 2.968  \\  \hline \hline
4 & 20    &  0.437--03& 0.358--03&   \\ 
  & 26    &  0.153--03& 0.126--03& 3.989  \\ 
  & 32    &  0.669--04& 0.551--04& 3.984  \\  \hline \hline
5 & 8     &  0.424--02& 0.346--02&   \\ 
  & 18    &  0.719--04& 0.583--04& 5.029  \\ 
  & 32    &  0.404--05& 0.327--05& 5.005  \\  \hline
\end{tabular}
}
\end{center}
\begin{center}
\textbf{Table 6.} Example 4: $H^1(\Omega)$ errors at $T=1$ for (\ref{eq:case2})
\end{center}
\end{table}
 \begin{table}[hbt]
\begin{center}
{\scriptsize
\begin{tabular}{|c|c|c|c|c|}
\hline $r$ & $N$  &$\|e_{h,1}\|_{L^\infty(\Omega)}$ & $\|e_{h,2}\|_{L^\infty(\Omega)}$ & Rate \\
\hline
3 & 20    &  0.175--03& 0.142--03&   \\ 
  & 26    &  0.621--04& 0.506--04& 3.951  \\ 
  & 32    &  0.272--04& 0.222--04& 3.973  \\  \hline \hline
4 &  8    &  0.106--02& 0.751--03&   \\ 
  & 18    &  0.176--04& 0.120--04& 5.055  \\ 
  & 32    &  0.963--06& 0.659--06& 5.048  \\  \hline \hline
5 & 20    &  0.102--05& 0.698--06&   \\ 
  & 26    &  0.214--06& 0.145--06& 5.977  \\ 
  & 32    &  0.616--07& 0.420--07& 5.984  \\  \hline
\end{tabular}
}
\end{center}
\begin{center}
\textbf{Table 7.} Example 4:  $L^\infty(\Omega)$ errors at $T=1$ for (\ref{eq:case2})
\end{center}
\end{table}
 \begin{table}[hbt]
\begin{center}
{\scriptsize
\begin{tabular}{|c|c|c|c|c|c|c|c|c|c|c|}
\hline
&    & \multicolumn{2}{c|}{Maximum nodal errors}  &   & \multicolumn{2}{c|}{Maximum nodal errors} &   & \multicolumn{2}{c|}{Maximum nodal errors} &  \\ \cline{3-4} \cline{6-7} \cline{9-10}
$r$ & $N$  & $e_{h,1}$ & $e_{h,2}$ & Rate & $\partial e_{h,1}/\partial x$ & $\partial e_{h,2}/\partial x$ & Rate & $\partial e_{h,1}/\partial y$ & $\partial e_{h,2}/\partial y$ & Rate\\
\hline
3 & 20 & 0.282--03& 0.284--03&       & 0.203--02& 0.894--03&      & 0.896--03& 0.204--02&  \\
  & 24 & 0.136--03& 0.137--03& 3.998 & 0.103--02& 0.432--03& 3.726& 0.432--03& 0.104--02& 3.727\\
  & 28 & 0.734--04& 0.739--04& 3.999 & 0.543--03& 0.233--03& 4.160& 0.233--03& 0.546--03& 4.159 \\
  & 32 & 0.430--04& 0.433--04& 3.999 & 0.327--03& 0.137--03& 3.812& 0.137--03& 0.328--03& 3.816\\
  & 36 & 0.269--04& 0.270--04& 3.999 & 0.201--03& 0.853--04& 4.126& 0.855--04& 0.202--03& 4.125 \\  \hline \hline
4 & 20 & 0.345--05& 0.346--05&       & 0.205--04& 0.108--04&      & 0.109--04& 0.206--04&      \\
  & 24 & 0.115--05& 0.116--05& 6.000 & 0.721--05& 0.363--05& 5.728& 0.364--05& 0.724--05& 5.729 \\
  & 28 & 0.458--06& 0.460--06& 6.000 & 0.279--05& 0.144--05& 6.162& 0.144--05& 0.280--05& 6.162 \\
  & 32 & 0.205--06& 0.206--06& 6.000 & 0.128--05& 0.647--06& 5.813& 0.647--06& 0.129--05& 5.813\\
  & 36 & 0.101--06& 0.102--06& 6.000 & 0.624--06& 0.319--06& 6.128& 0.319--06& 0.626--06& 6.127 \\  \hline \hline
5 & 20 & 0.340--07& 0.342--07&       & 0.202--06& 0.107--06&      & 0.107--06& 0.203--06&  \\
  & 24 & 0.792--08& 0.795--08& 8.000 & 0.495--07& 0.248--07& 7.728& 0.249--07& 0.497--07& 7.729\\
  & 28 & 0.231--08& 0.232--08& 8.000 & 0.141--07& 0.724--08& 8.162& 0.724--08& 0.141--07& 8.161 \\
  & 32 & 0.792--09& 0.796--09& 8.000 & 0.495--08& 0.247--08& 7.813& 0.249--08& 0.497--08& 7.813\\
  & 36 & 0.313--09& 0.315--09& 7.873 & 0.193--08& 0.984--09& 8.000& 0.985--09& 0.194--08& 8.000 \\
\hline
 \end{tabular}
 }
\end{center}
\begin{center}
\textbf{Table 8.} Example 4: Maximum nodal errors at $T=1$ for (\ref{eq:case1})
\end{center}
\end{table}


In Figure 5, we present graphs of $u_{h,1}$ and $u_{h,2}$ at $T=1$ computed with $N=20$, (that is, $h=0.1$), $\tau=h^3$ and $r=4$. When compared with corresponding graphs (Figures  5a, 6a and  Figures 5b, 6b) in \cite{ZhWoZh}, we see that the graphs are similar, taking into account their differing orientations.
Figure 6 shows graphs of the corresponding nodal errors $e_{h,1}$ and $e_{h,2}$ at $T=1$.
It should be noted that the nodal errors in this figure
are approximately $10^{-6}$
which are of the same order as those obtained in \cite[Figures 5d, 6d]{ZhWoZh}.

\begin{center}
\resizebox{1\linewidth}{0.3\linewidth}{\includegraphics{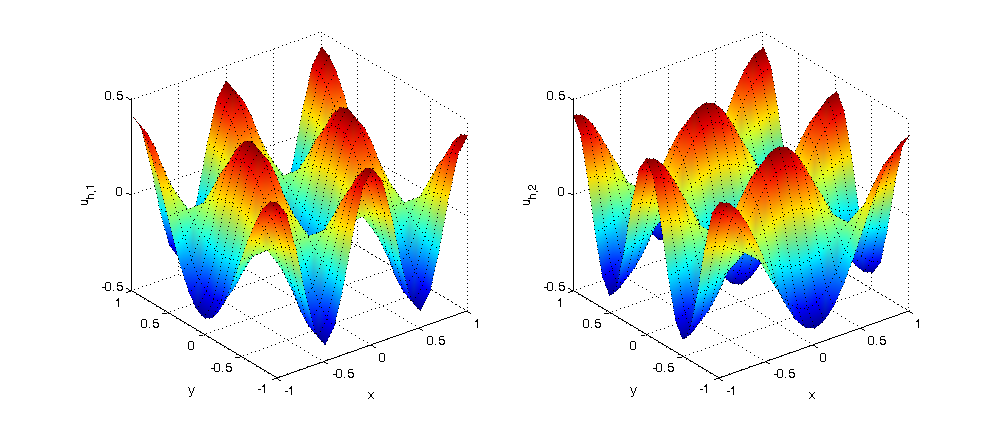}}
\end{center}
\begin{center}
{\bf Figure 5.} Example 4: The approximate solutions $u_{h,1}$ and $u_{h,2}$ at $T=1$
\end{center}

\begin{center}
\resizebox{1\linewidth}{0.4\linewidth}{\includegraphics{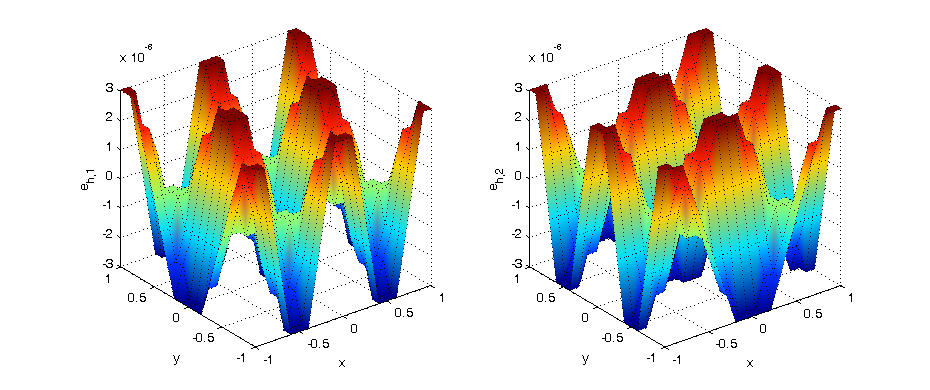}}
\end{center}
\begin{center}
{\bf Figure 6.}  Example 4: The errors  $e_{h,1}$ and $e_{h,2}$ at $T=1$
\end{center}

\subsection{Gierer-Meinhardt Model}
The Gierer-Meinhardt model \cite{GiMe} comprises
(\ref{prob1})--(\ref{prob4}) with
\beq
\label{qrhs}
{\bf f}({\bf u}) = [u_1^2/u_2 - u_1, u_1^2/(\epsilon\mu) -u_2/\mu]^T.
\eeq
First we examine the performance of the ADI  method on test problems considered in  \cite{McDoGi}
in which $\Omega = (-1,1)\times (-1,1)$, $\mu > 0$, and $D_1 = \epsilon^2$, $D_2= \kappa/\mu$,
where $\epsilon = 0.04$, $\mu=0.1$, and $\kappa$ is varied.
Also, ${\bf u}(x,y,0) = [g^0_1(x,y),g^0_2(x,y)]^T$, where%
\beq
\label{qic}
\dsp{ g_1^0=\frac{1}{2} \left[1+0.001 \sum_{k=1}^{20} \cos \left( \frac{k \pi y}{2} \right)\right]
\mbox{sech}^2\left(\frac{\sqrt{x^2+y^2}}{2 \epsilon}\right),}
\qquad
\dsp{g_2^0=\frac{\cosh\left(1-\sqrt{x^2+y^2}\right)}{3 \cosh(1)}}.
\eeq
In \cite{McDoGi}, this problem, which has no known closed-form solution, is solved using a Chebyshev spectral collocation method for the spatial discretization and a linearized backward Euler method for the time--stepping. At each time step,  the collocation equations are solved using a preconditioned GMRES method at a cost of $O({\cal N}^{3/2})$ operations, where $\cal N$ is the number of unknowns.  As in \cite{McDoGi}, we restrict our attention to the dynamics of $u_1$, the initial profile
of which is shown in Figure~7.
\begin{center}
\resizebox{0.6\linewidth}{0.35\linewidth}{\includegraphics{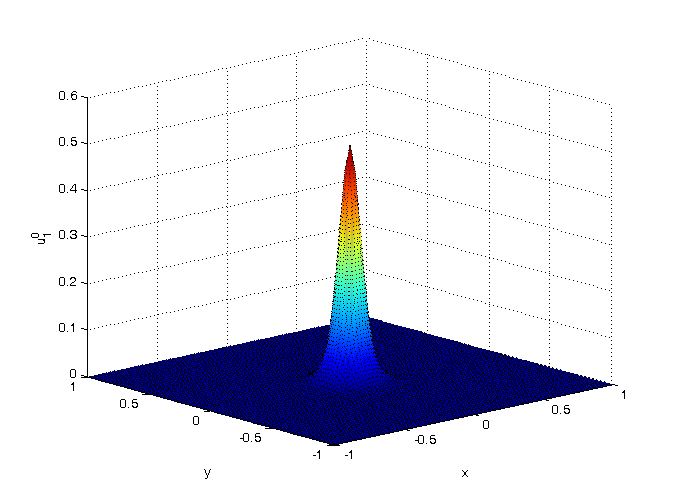}}
\end{center}
\begin{center}
{\bf Figure 7.} Initial profile of  $u_1^0$
\end{center}

\vspace{.2in}
\noindent
{\bf Example 5.}
In this example, $\kappa=0.0128$.
For $r=3$, $N=20$ (that is, $h=0.1$), $\tau=h^2$, we graph  $u_{h,1}$ at the values of $t$ selected in \cite{McDoGi}. Figure 8 provides graphs corresponding to Figures 7--10 in \cite{McDoGi}.
It is clear from our graphs that the dynamics of $u_{h,1}$ follow the same general pattern displayed in \cite{McDoGi}. The center of the spike at the origin sinks to the floor producing a ring formation that moves outward to the boundary ($t=320$). Then it collapses into smaller spikes at the boundary and generates four new spikes in the interior near the center of each of the four boundaries ($t=420$). Thereafter it appears that the spikes evolve to a symmetric pattern ($t=900$). This pattern is more obvious in Figure 9 which shows aerial views of the graphs in Figure 8. There are, however, some inexplicable differences between our graphs and those in \cite{McDoGi}. For instance, if one compares corresponding graphs at $t=900$, our graph shows an arrangement with four spikes along each boundary and four in the middle whereas in \cite{McDoGi} there appear to be five spikes along the $y$-boundaries, four along the $x$-boundaries and four in the middle. Also, at $t=320, 340$, the spikes near the boundary appear to be more prominent than in our corresponding graphs.

\vspace{.2in}
\noindent
{\bf Example 6.}
In this example,  $\kappa=0.0152$ and we compare results obtained using the ADI method with the graphs in Figures 11--14 of \cite{McDoGi}. In our computations, we use the same parameters as in Example 5, with $h=0.1$, $\tau=h^2$ and $r=3$. Figure 10 shows graphs of $u_{h,1}$ for various values of $t$. Once again, in our graphs, we observe the general pattern obtained for $u_{h,1}$ in Figures 11--14 of \cite{McDoGi}. The spike splits into two spikes spreading in the $x$ direction and becomes symmetric ($t=140$). Then, each of the spikes splits, spreading in the $y$ direction, and maintains symmetry ($t=290$). Next, each of the four spikes splits into two along the $x$ direction, and the eight spikes arrange themselves symmetrically about the center ($t=620$). Finally, each of the four outermost spikes split into two and the 12 spikes arrange themselves symmetrically about the center ending with four spikes at the corners and eight spikes in a circular pattern around the center ($t=990$). This is seen more clearly in Figure 11 which gives aerial views of the graphs of Figure 10. Comparing our pattern with that of \cite{McDoGi}, we see a similarity up to $t=520$.
At $t=570$, it appears that in \cite{McDoGi} the inner spikes split as opposed to the outer ones in our case. Thus, while at $t=990$ there are 12 spikes in our graph and in \cite{McDoGi},  the symmetric pattern that we obtain is not seen in the corresponding graph of \cite{McDoGi}.

\begin{center}
\resizebox{1\linewidth}{1.3\linewidth}{\includegraphics{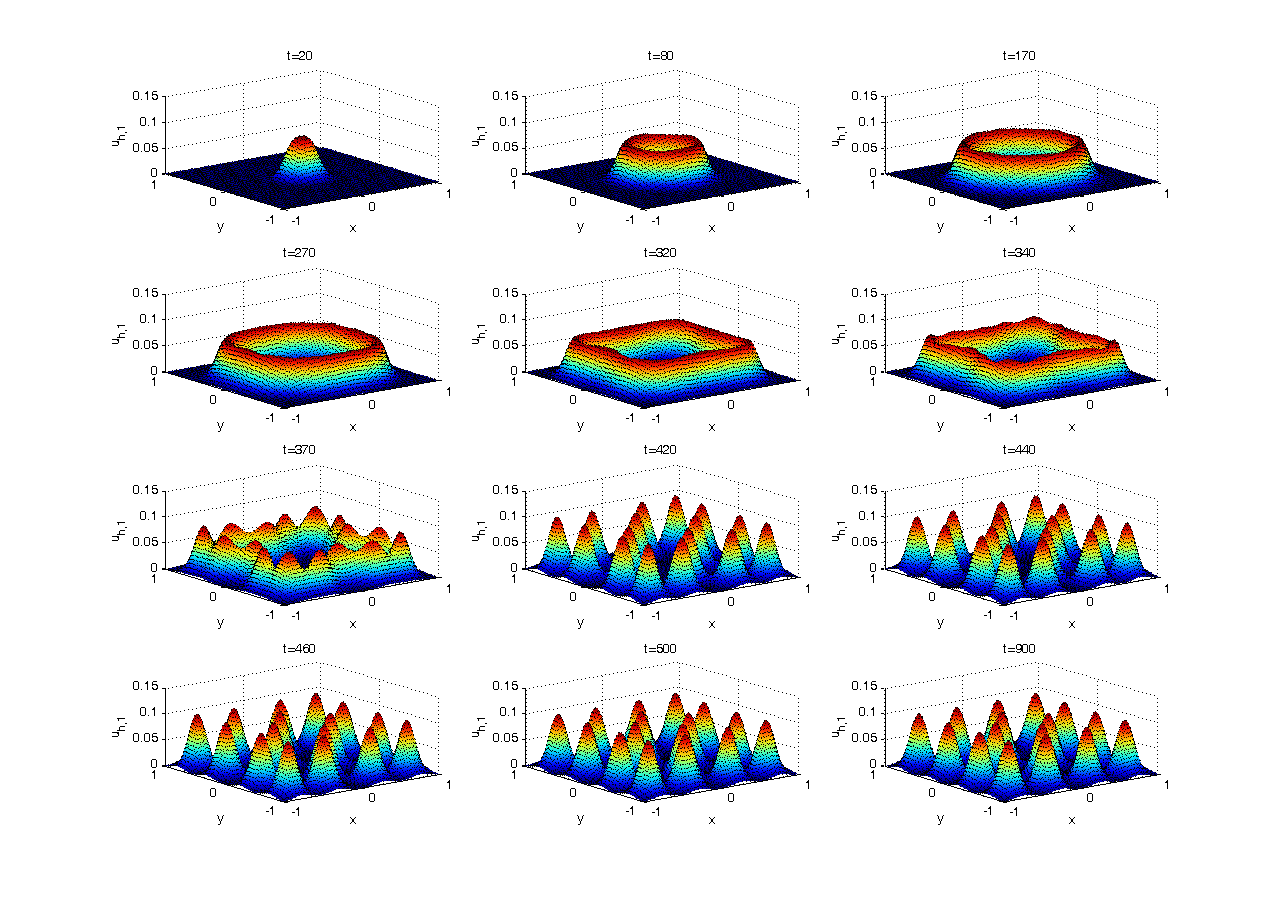}}

{\bf Figure 8.} Example 5: Graphs of $u_{h,1}$ at various values of $t$ 
\resizebox{1\linewidth}{1\linewidth}{\includegraphics{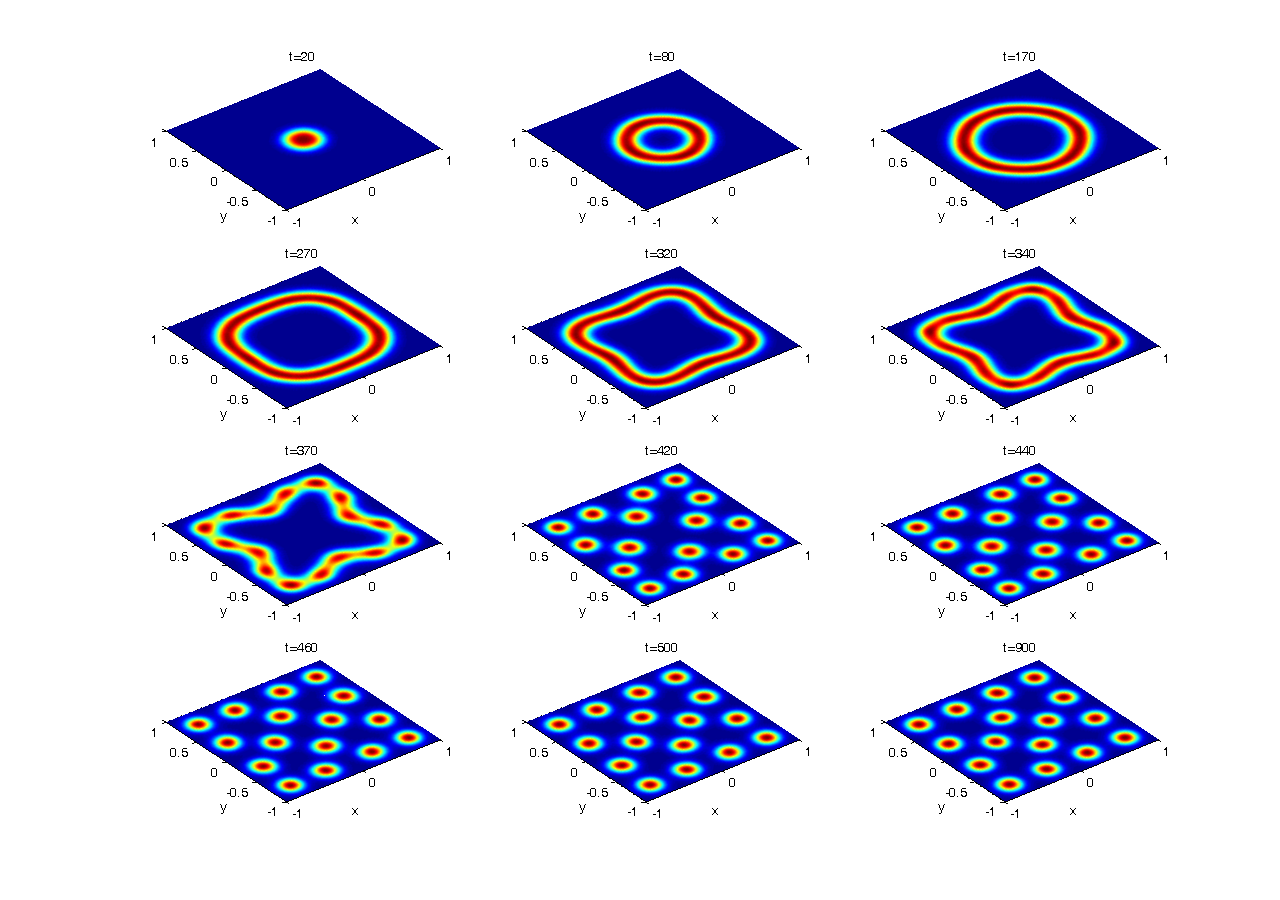}}

{\bf Figure 9.} Example 5: Aerial views of $u_{h,1}$ at various values of $t$
\end{center}

\begin{center}
\resizebox{1\linewidth}{1.1\linewidth}{\includegraphics{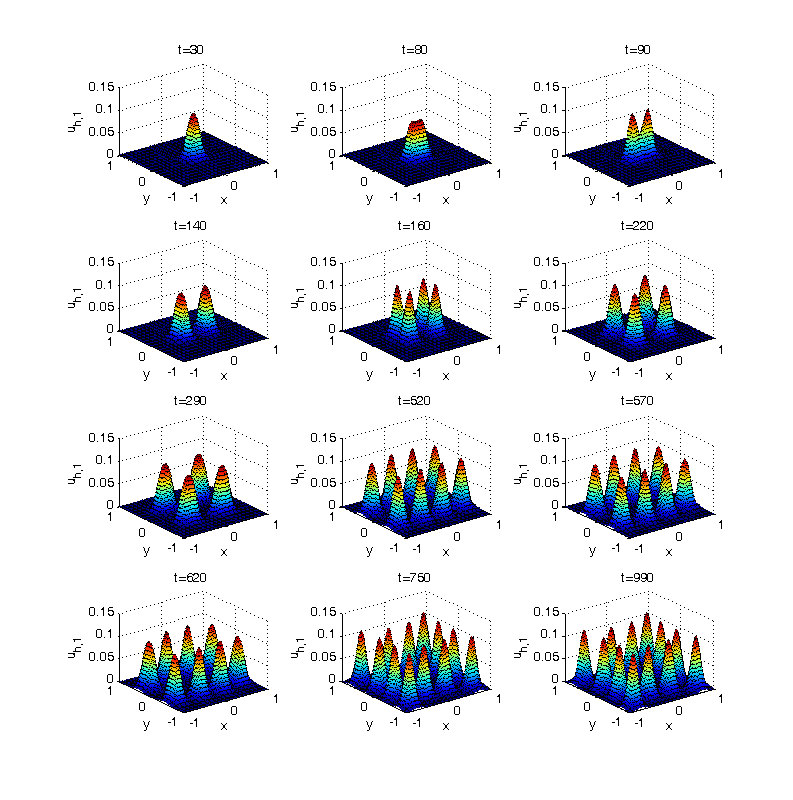}}

{\bf Figure 10.} Example 6: Graphs of $u_{h,1}$ at various values of $t$

\resizebox{1\linewidth}{1\linewidth}{\includegraphics{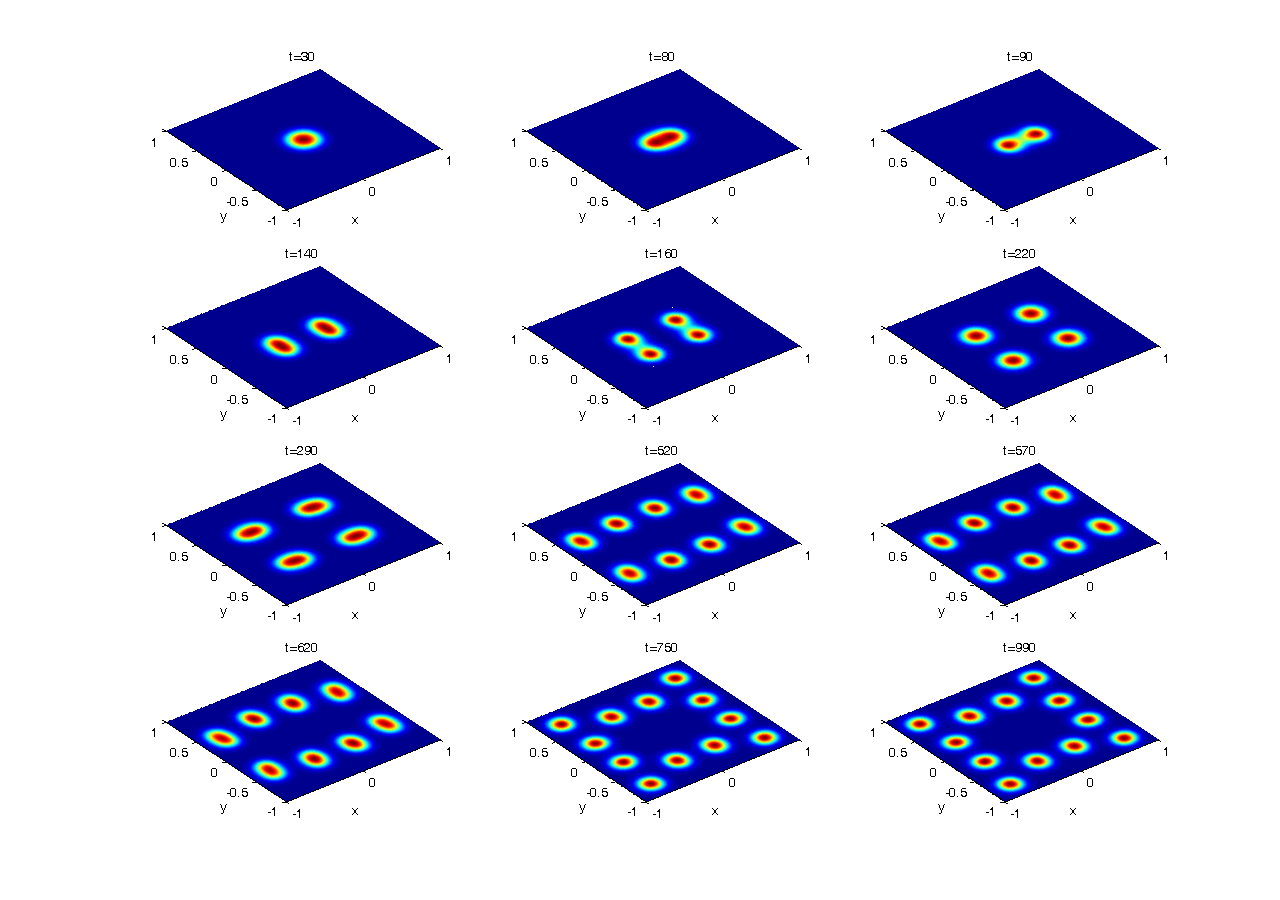}}
 \end{center}
 \begin{center}

{\bf Figure 11.} Example 6: Aerial views of $u_{h,1}$ at various values of $t$
\end{center}
On comparing Examples 5 and 6, we see that a modest change in the value of $\kappa$ leads to a significant change in the behavior of the solution.

Examples 5 and 6, taken from \cite{McDoGi}, are based on examples in \cite[Section 4.2]{Qi}. While the solution in \cite{McDoGi} is a simple scaling of the solution in  \cite[(4.2)-(4.6)]{Qi}, the initial value prescribed is not the scaled initial value used in \cite{Qi}; each of the components in the initial condition should be divided by $\epsilon$. Moreover, the graphs obtained in \cite{McDoGi} are not given at the same $t$ values as in \cite{Qi}. In view of these discrepancies, the comparison given in \cite{McDoGi} between the results presented therein and those of \cite{Qi} is questionable.
Consequently, we next consider the examples presented in \cite[Section 4.2]{Qi}, where the technique employed is based on a moving mesh finite element Galerkin method with piecewise linear functions on triangulations of $\Omega$, with a second-order Runge-Kutta scheme for the time-stepping. The test problem in this case is ({\ref{prob1}) with {\bf f}$(${\bf u}$)$ given by (\ref{qrhs}) in which $\epsilon$ is replaced by $\epsilon^2$. The initial value (which is incorrectly stated in \cite[(4.1)]{Qi}; $\epsilon^2$ should be $\epsilon$) is given by (\ref{qic}). All other parameter choices are the same as in Examples 5 and 6.
\begin{center}
\resizebox{1\linewidth}{0.6\linewidth}{\includegraphics{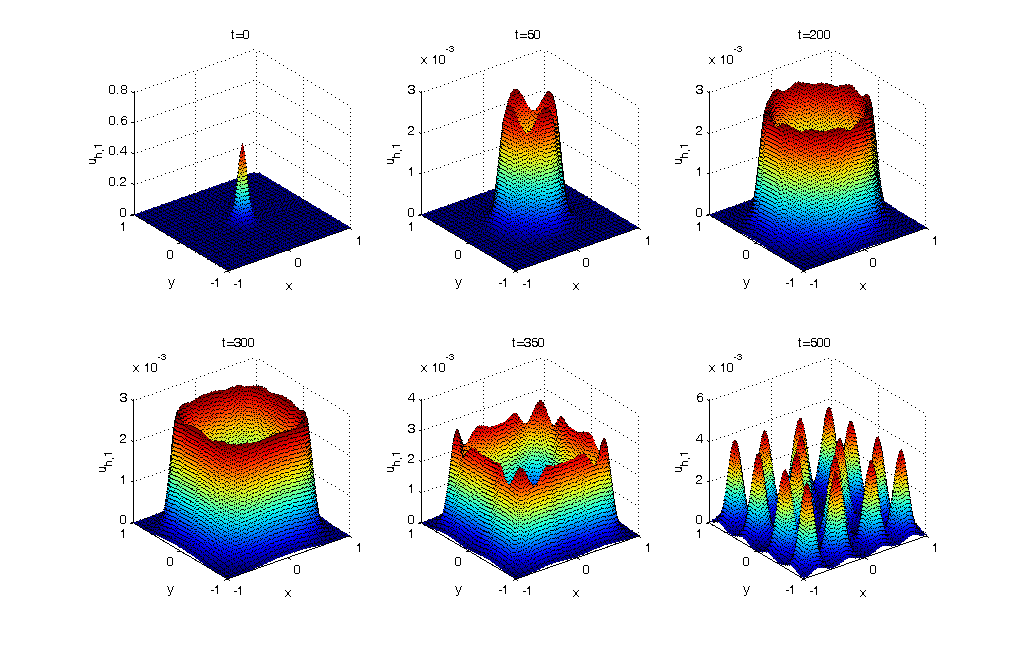}}
\end{center}
\begin{center}
{\bf Figure 12.} Example 7: Graphs of $u_{h,1}$ at various values of $t$
\end{center}

\begin{center}
\resizebox{1\linewidth}{0.5\linewidth}{\includegraphics{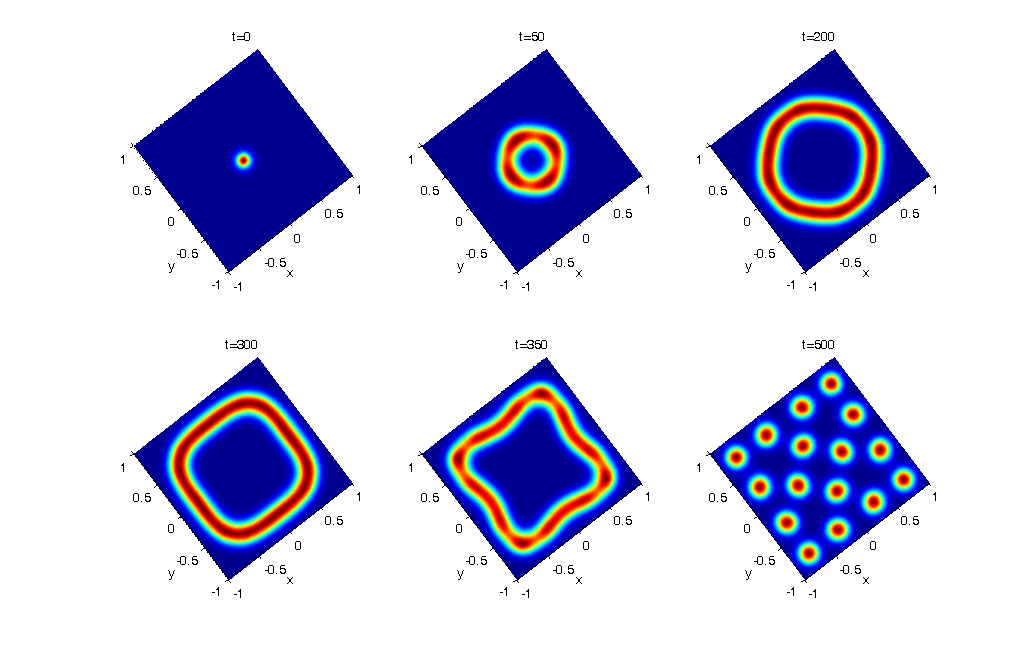}}
\end{center}
\begin{center}
{\bf Figure 13.} Example 7: Aerial views of $u_{h,1}$ at various values of $t$
\end{center}

\vspace{.2in}
\noindent
{\bf Example 7.}
In this example, $\kappa=0.0128$.
For $r=3$, $N=20$ (that is, $h=0.1$), $\tau=h^2$, we graph  $u_{h,1}$ at the values of $t$ selected in \cite{Qi}. Figures 12 and 13 show the graphs and aerial views, respectively, corresponding to Figures 9--10 in \cite{Qi}. There is fairly good agreement between the two sets of graphs.

\vspace{.2in}
\noindent
{\bf Example 8.}
In this example,  $\kappa=0.0152$ and we compare results obtained using the ADI method with the graphs in Figures 11--12 of \cite{Qi}. 
With $h=0.1$, $\tau=h^2$ and $r=3$, Figure 14 shows graphs of $u_{h,1}$ for various values of $t$. From the aerial view of the solution shown in Figure 15, it appears that we have good agreement up to approximately $t=100$. Thereafter, there is a difference in the arrangement of the spikes. Our figures indicate that a steady-state pattern is reached by $t=500$ which is not the case in \cite{Qi}.

Again, we see that a small change in $\kappa$ leads to a significant change in the behavior of the solution.

\begin{center}
\resizebox{1\linewidth}{0.8\linewidth}{\includegraphics{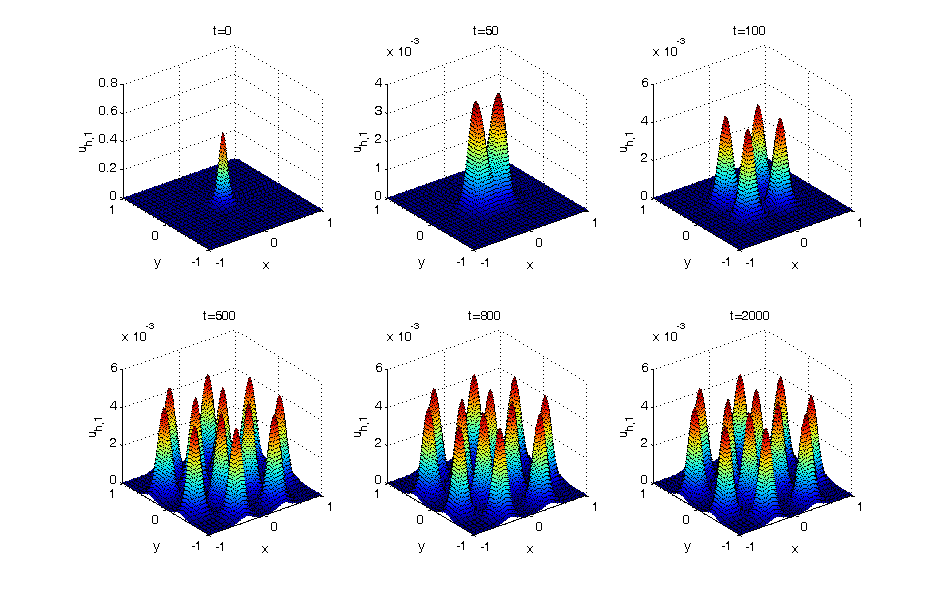}}
\end{center}
\begin{center}
{\bf Figure 14.} Example 8: Graphs of $u_{h,1}$ at various values of $t$
\end{center}

\begin{center}
\resizebox{1\linewidth}{0.5\linewidth}{\includegraphics{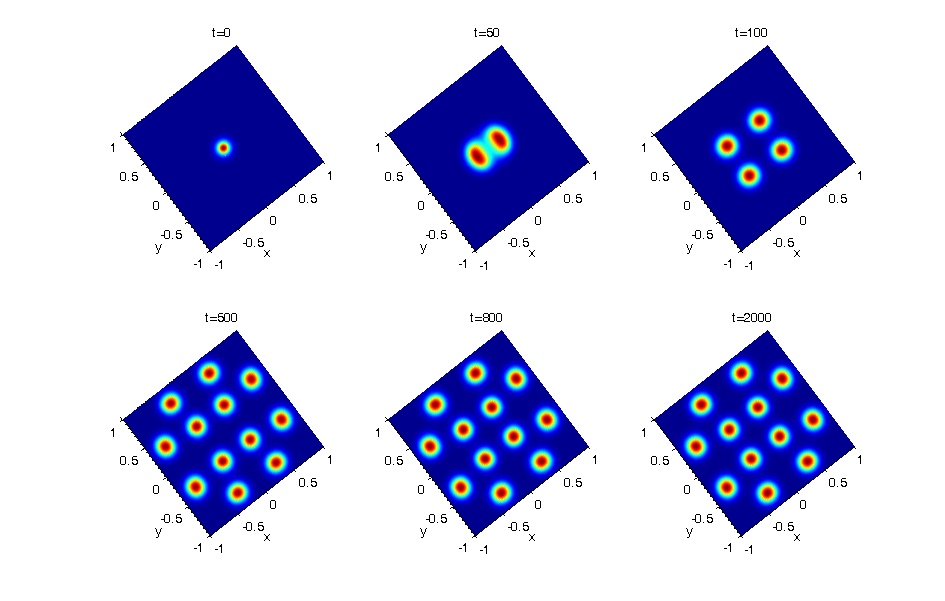}}
\end{center}
\begin{center}
{\bf Figure 15.} Example 8: Aerial views of $u_{h,1}$ at various values of $t$
\end{center}

\subsection{Schnakenberg Model}
\label{sm}
The Schnakenberg model \cite{Sc} 
is given by (\ref{prob1})--(\ref{prob4})
with
\[
{\bf f}({\bf u}) = [\gamma(a - u_1 +u_1^2u_2), \gamma(b-u_1^2u_2)]^T,
\]
where $\gamma, a$ and $b$ are constants.
In \cite{Ma}, this model is studied on fixed and growing domains using methods based on
a piecewise linear finite element Galerkin discretization in space coupled with a linearized backward Euler method or a linearized second-order BDF for the time discretization. The algebraic systems arising at each time step are solved iteratively using a preconditioned GMRES method.

In the following test problems, we consider the fixed domain
$
\Omega = (0,1) \times (0,1),
$
and take
$
D_1=1,  \quad D_2=10.
$
 \begin{table}[hbt]
\begin{center}
{\scriptsize
\begin{tabular}{|c|c|c|c|c|}
\hline $r$ & $N_l$  & $\|e_{h,1}\|_{L^2(\Omega)}$ & $\|e_{h,2}\|_{L^2(\Omega)}$ & Rate \\
\hline
3 & 10    &  0.683--04& 0.805--03&   \\ 
  & 15    &  0.134--04& 0.159--03& 3.998  \\ 
  & 20    &  0.424--05& 0.504--04& 4.000  \\  \hline \hline
4 &  9    &  0.132--04& 0.141--03&   \\ 
  & 16    &  0.746--06& 0.792--05& 4.999  \\ 
  & 25    &  0.802--07& 0.851--06& 4.999  \\  \hline \hline
5 & 10    &  0.775--06& 0.831--05&   \\ 
  & 15    &  0.680--07& 0.730--06& 6.000  \\ 
  & 20    &  0.121--07& 0.130--06& 6.000  \\  \hline
\end{tabular}
}
\end{center}
\begin{center}
\textbf{Table 10.} Example 9: $L^2(\Omega)$ errors at $T=1$
\end{center}
\end{table}
 \begin{table}[hbt]
\begin{center}
{\scriptsize
\begin{tabular}{|c|c|c|c|c|}
\hline $r$ & $N_l$  &$\|e_{h,1}\|_{H^1(\Omega)}$ & $\|e_{h,2}\|_{H^1(\Omega)}$& Rate \\
\hline
3 & 9     &  0.764--02& 0.792--01&   \\ 
  & 16    &  0.138--02& 0.143--01& 2.971  \\ 
  & 25    &  0.363--03& 0.377--02& 2.993  \\  \hline \hline
4 & 10    &  0.527--03& 0.589--02&   \\ 
  & 15    &  0.104--03& 0.116--02& 3.997  \\ 
  & 20    &  0.330--04& 0.369--03& 3.999  \\  \hline \hline
5 & 9     &  0.882--04& 0.999--03&   \\ 
  & 16    &  0.497--05& 0.563--04& 5.000  \\ 
  & 25    &  0.533--06& 0.604--05& 5.000  \\  \hline
\end{tabular}
}
\end{center}
\begin{center}
\textbf{Table 11.} Example 9: $H^1(\Omega)$ errors  at $T=1$
\end{center}
\end{table}
 \begin{table}[hbt]
\begin{center}
{\scriptsize
\begin{tabular}{|c|c|c|c|c|}
\hline $r$ & $N_l$  & $\|e_{h,1}\|_{L^{\infty}(\Omega)}$ & $\|e_{h,2}\|_{L^{\infty}(\Omega)}$ & Rate \\
\hline
3 & 10    &  0.306--03& 0.174--02&   \\ 
  & 15    &  0.604--04& 0.341--03& 4.021  \\ 
  & 20    &  0.191--04& 0.109--03& 3.968  \\  \hline \hline
4 &  9    &  0.393--04& 0.283--04&   \\ 
  & 16    &  0.223--05& 0.159--04& 5.003  \\ 
  & 25    &  0.241--06& 0.171--05& 4.998  \\  \hline \hline
5 & 10    &  0.238--05& 0.167--04&   \\ 
  & 15    &  0.209--06& 0.146--05& 6.009  \\ 
  & 20    &  0.373--07& 0.262--06& 5.982  \\  \hline
\end{tabular}
}
\end{center}
\begin{center}
\textbf{Table 12.} Example 9: $L^\infty(\Omega)$ errors at $T=1$
\end{center}
\end{table}

 \begin{table}[hbt]
\begin{center}
{\scriptsize
\begin{tabular}{|c|c|c|c|c|c|c|c|c|c|c|}
\hline
&    & \multicolumn{2}{c|}{Maximum nodal errors}  &   & \multicolumn{2}{c|}{Maximum nodal errors} &   & \multicolumn{2}{c|}{Maximum nodal errors} &  \\ \cline{3-4} \cline{6-7} \cline{9-10}
$r$ & $N_l$  & $e_{h,1}$ & $e_{h,2}$ & Rate & $\partial e_{h,1}/\partial x$ & $\partial e_{h,2}/\partial x$ & Rate & $\partial e_{h,1}/\partial y$ & $\partial e_{h,2}/\partial y$ & Rate\\
\hline
3 & 10 & 0.306--03& 0.174--02&       & 0.106--02& 0.546--02&      & 0.699--03& 0.989--02&  \\
  & 15 & 0.604--04& 0.341--03& 4.021 & 0.209--03& 0.106--02& 4.033& 0.137--03& 0.204--02& 3.888\\
  & 20 & 0.191--04& 0.109--03& 3.968 & 0.661--04& 0.341--03& 3.952& 0.438--04& 0.650--03& 3.981 \\  \hline \hline
4 & 10 & 0.227--05& 0.166--04&       & 0.104--04& 0.522--04&      & 0.519--05& 0.992--04&  \\
  & 15 & 0.199--06& 0.145--05& 6.024 & 0.913--06& 0.452--05& 6.035& 0.458--06& 0.911--05& 5.890\\
  & 20 & 0.354--07& 0.260--06& 5.965 & 0.163--06& 0.816--06& 5.950& 0.823--07& 0.163--05& 5.981 \\  \hline \hline
5 & 10 & 0.233--07& 0.167--06&       & 0.108--06& 0.525--06&      & 0.529--07& 0.994--06&  \\
  & 15 & 0.909--09& 0.645--08& 8.025 & 0.420--08& 0.202--07& 8.036& 0.209--08& 0.405--07& 7.891\\
  & 20 & 0.924--10& 0.654--09& 7.958 & 0.429--09& 0.205--08& 7.943& 0.213--09& 0.409--08& 7.974 \\
\hline
 \end{tabular}
 }
\end{center}
\begin{center}
\textbf{Table 13.} Example 9: Maximum nodal errors  at $T=1$
\end{center}
\end{table}

\noindent
{\bf Example 9.}
In this example, we set $\gamma=10$, $a=0.1$, $b=0.9$, and construct
the reaction kinetic functions and initial functions so that the exact solution is given by (\ref{exsol}). Choosing the time step for various values of $r$ as
in subsection \ref{bm}, we present in Tables 10--13 the $H^k(\Omega), k=0,1,$ and $L^\infty(\Omega)$ norms of the errors $e_{h,1}$ and $e_{h,2}$ and maximum nodal errors at $T=1$, respectively.
These numerical results confirm the expected optimal rates of convergence and superconvergence at the nodes.

\vspace{.2in}
\noindent
{\bf Example 10.}
This test problem is Example 4 of \cite{Ma}, (cf.,  \cite[Example 5]{Ru}) in which
$\gamma=1000, a=0.126779, b=0.792366$, and
\beas
 g_1^0=0.919145+0.0016 \cos(2 \pi(x+y)) +0.01 \sum_{j=1}^8 \cos(2\pi j x),\\
g_2^0=0.937903+0.0016 \cos(2 \pi(x+y)) +0.01 \sum_{j=1}^8 \cos(2\pi j x).
\eeas
In this example, we take $r=5, N=20$ (that is, $h=0.1$),  and  $\tau=h^{3}$.
\begin{center}
\resizebox{1\linewidth}{1.2\linewidth}{\includegraphics{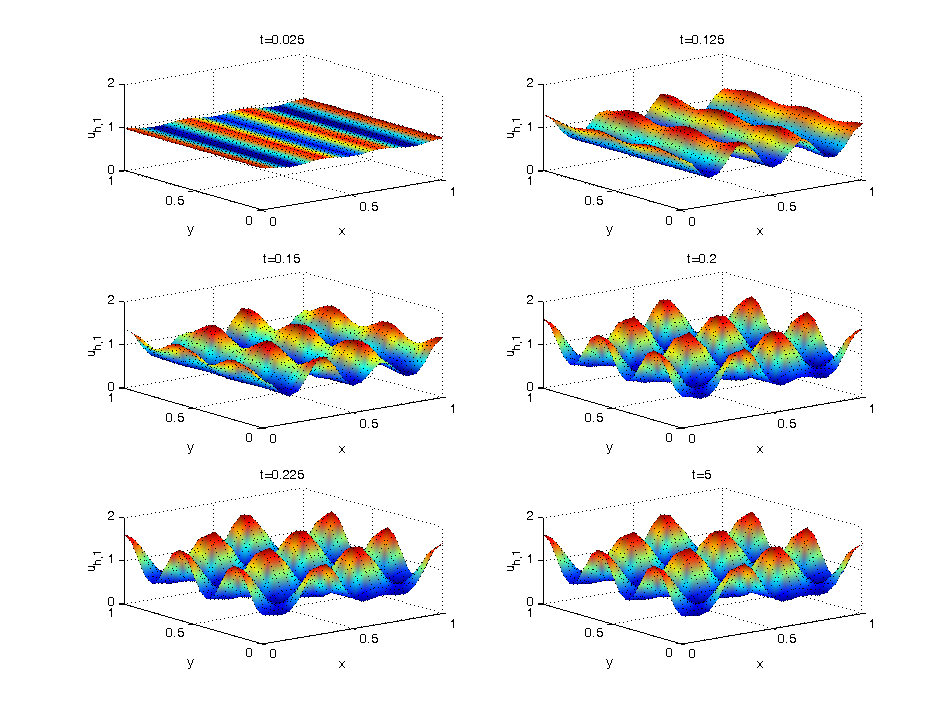}}
\end{center}
\begin{center}
{\bf Figure 16.} Example 10: Graphs of $u_{h,1}$ at various values of $t$
\end{center}
\begin{center}
\resizebox{1\linewidth}{0.9\linewidth}{\includegraphics{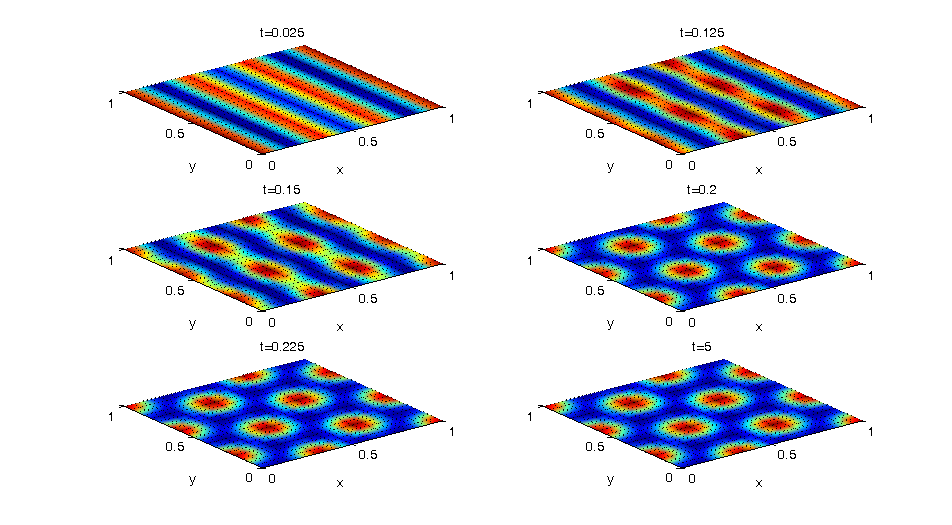}}
\end{center}
\begin{center}
{\bf Figure 17.} Example 10: Aerial views of $u_{h,1}$ at various values of $t$
\end{center}

Graphs of $u_{h,1}$ at various values of $t$ are presented in Figure 16 and corresponding aerial views are shown in Figure 17. On comparing the graphs of Figure 17 (rotated clockwise) with those of Figure 8 or 9 of \cite{Ma}, it is clear that there is agreement between them when $t \geq 0.2$. For $t < 0.2$, a similarity is difficult to detect
as the details in the graphs of \cite{Ma} are obscured because of the use of constant threshold shading.

\vspace{.2in}
\noindent
{\bf Example 11.}
Here we consider Example 5 of \cite{Ma}
where $\gamma=10000, a=-0.887757, b=2.774242$, and
\[
g_1^0=1.886485+0.001 \sum_{j=1}^{37} \frac{\cos(2\pi j x)}{j},\qquad
g_2^0=0.779539+0.001 \sum_{j=1}^{37} \frac{\cos(2\pi j x)}{j}.
\]
In this example, we take $r=6, N=20$ (that is, $h=0.05$), and $\tau=100 h^{3}$ so that the number of time steps, $N_t$, for the various values of $t$ selected in \cite{Ma} is an integer.

Graphs of $u_{h,1}$ at various values of $t$ are presented in Figure 18 and corresponding aerial views are shown in Figure 19. On comparing the graphs of Figure 19 (rotated clockwise) with those of Figure 13 in \cite{Ma}, we do not observe the ripples obtained in \cite{Ma}.

\begin{center}
\resizebox{1\linewidth}{!}{\includegraphics{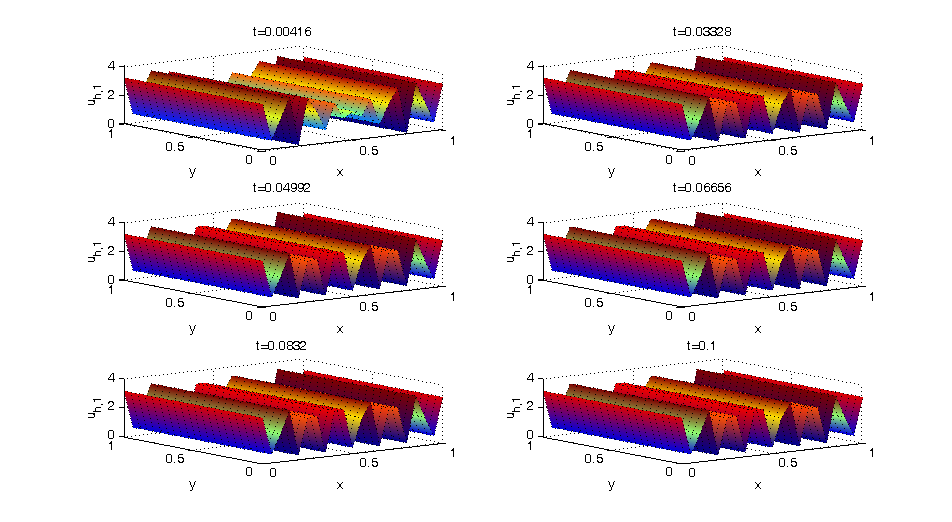}}
\end{center}
\begin{center}
{\bf Figure 18.} Example 11: Graphs of $u_{h,1}$ at various values of $t$
\end{center}
\begin{center}
\resizebox{1\linewidth}{!}{\includegraphics{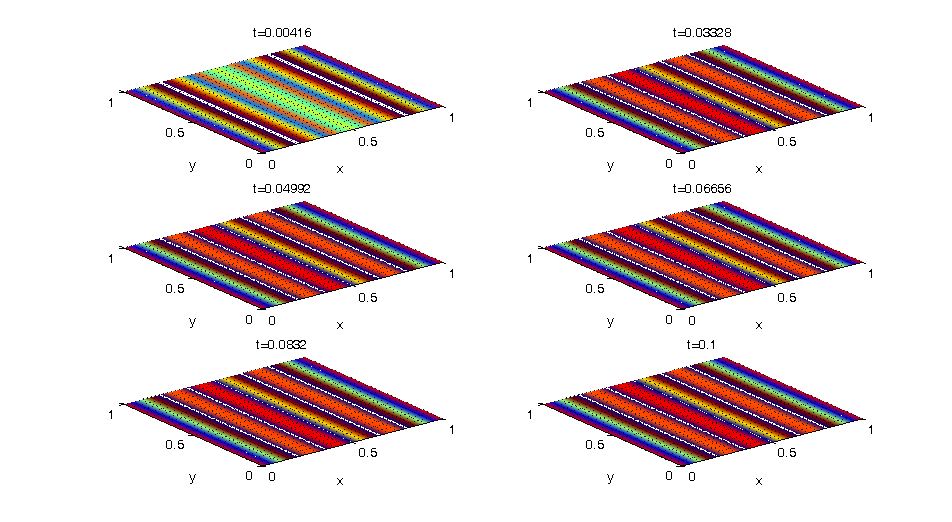}}
\end{center}
\begin{center}
{\bf Figure 19.} Example 11: Aerial views of $u_{h,1}$ at various values of $t$
\end{center}

\section{Concluding Remarks}
\label{sec:5}
\setcounter{equation}{0}
We have presented an  ADI  method for solving nonlinear systems of reaction-diffusion problems of the form (\ref{prob1}) which is second-order accurate in time and of optimal accuracy in space. The method enjoys several features that render it more efficient than current methods for solving such problems. Using test problems primarily from the literature, we have examined the accuracy of the method and demonstrated optimal convergence rates in standard norms for the Brusselator, Gray-Scott and Schnakenberg models. Also, we have compared results obtained by the ADI scheme with results in the literature for the Brusselator, Gierer-Meinhardt and Schnakenberg models for which exact solutions are not known. In the process, we have identified several errors and discrepancies in the literature.

The ADI  method can be generalized to solve a system of nonlinear parabolic equations with more than two equations.
It is also possible to generalize the scheme to systems of nonlinear parabolic problems on rectangles where the diffusion coefficients and the 
reaction kinetic functions depend on $x,y,t,{\bf u},\nabla {\bf u}$ as well as when the boundary conditions are of Dirichlet, Neumann or mixed type; cf., \cite{BiFe1}. In view of the high accuracy that can be achieved in space, future research will include the development of time-stepping methods of higher order accuracy than second; cf., \cite{BiFe2}. Also, the treatment of growing domains (cf., \cite{Ma}) and the formulation of ADI OSC methods for problems in more general regions will be considered; cf., \cite{An,Le,LiTa}.

\section*{Acknowledgement}
The authors wish to thank Professor Andreas Karageorghis of the University of Cyprus for his assistance during the preparation of this paper.

\end{document}